\documentclass[aos,MSNbibl,nameyear,seceqn,dvips]{arximspdf}
\usepackage{graphicx}

\doi{10.1214/11-AOS890}
\volume{39}
\issue{4}
\pubyear{2011}
\firstpage{1963}
\lastpage{2006}

\makeatletter

\newtheorem{theorem}{Theorem}[section]
\newtheorem{lemma}[theorem]{Lemma}

\newproclaim{example}[theorem]{Example}
\newproclaim{remark}[theorem]{Remark}
\newtheorem{proposition}[theorem]{Proposition}

\newcommand{\Nat}{\mathbb{N}}
\newcommand{\Dm}{\mathbb{D}}
\newcommand{\PP}{\mathbb{P}}

\newcommand{\ind}{\mathbb{I}}
\newcommand{\B}{\mathbb{B}}
\newcommand{\A}{\mathbb{A}}
\newcommand{\G}{\mathbb{G}}
\newcommand{\Exp}{\mathbb{E}}

\newcommand{\x}{\mathbf{x}}
\newcommand{\y}{\mathbf{y}}
\newcommand{\weak}{\stackrel{w}{\rightsquigarrow}}
\newcommand{\weakP}{\mathop{\rightsquigarrow}_{\xi}^{\PP}}
\newcommand{\eps}{\varepsilon}
\newcommand{\W}{\mathbb{W}}

\newcommand{\weakPi}[1]{\mathop{\rightsquigarrow}_{#1}^{\PP}}
\newcommand{\R}{\mathbb{R}}
\newcommand{\Var}{\operatorname{Var}}
\newcommand{\Cov}{\operatorname{Cov}}

\makeatother

\begin{document}
\begin{frontmatter}

\title{New estimators of the Pickands dependence function and a test for extreme-value dependence\thanksref{T1}}
\runtitle{New estimators of Pickands dependence function}

\thankstext{T1}{Supported in part by the Collaborative Research Center ``Statistical modeling of nonlinear dynamic
processes'' (SFB 823) of the German Research Foundation (DFG).}

\begin{aug}
\author[A]{\fnms{Axel} \snm{B\"{u}cher}\ead[label=e1]{axel.buecher@ruhr-uni-bochum.de}},
\author[A]{\fnms{Holger} \snm{Dette}\corref{}\ead[label=e2]{holger.dette@ruhr-uni-bochum.de}}
and
\author[A]{\fnms{Stanislav} \snm{Volgushev}\ead[label=e3]{stanislav.volgushev@ruhr-uni-bochum.de}}

\runauthor{A. B\"{u}cher, H. Dette and S. Volgushev}
\affiliation{Ruhr-Universit\"{a}t Bochum}
\address[A]{Fakult\"{a}t f\"{u}r Mathematik\\
Ruhr-Universit\"{a}t Bochum\\
Universit\"{a}tsstra\ss e 150\\
44780 Bochum\\
Germany\\
\printead{e1}\\
\phantom{E-mail: }\printead*{e2}\\
\phantom{E-mail: }\printead*{e3}} 
\end{aug}

\received{\smonth{12} \syear{2010}}

%
\begin{abstract}
We propose a new class of estimators for Pickands dependence function
which is based on the concept of minimum distance estimation. An explicit
integral representation of the function $A^*(t)$, which minimizes a weighted
$L^2$-distance between the logarithm of the copula $C(y^{1-t},y^t)$ and
functions of the form $A(t) \log(y)$ is derived. If the unknown copula
is an
extreme-value copula, the function $A^*(t)$ coincides with Pickands
dependence function. Moreover, even if this is not the case, the function
$A^*(t)$ always satisfies the boundary conditions of a Pickands dependence
function. The estimators are obtained by replacing the unknown copula
by its
empirical counterpart and weak convergence of the corresponding process
is shown. A comparison with the commonly used estimators is performed
from a theoretical point of view and by means of a simulation study. Our
asymptotic and numerical results indicate that some of the new estimators
outperform the estimators, which were recently proposed by Genest and Segers
[\textit{Ann. Statist.} \textbf{37} (2009) 2990--3022].
As a by-product of our results, we obtain a simple test for the
hypothesis of
an extreme-value copula, which is consistent against all positive quadrant
dependent alternatives satisfying weak differentiability assumptions of
first \mbox{order}.
\end{abstract}

%

\end{frontmatter}

\section{Introduction}\label{sec1}

The copula provides
an elegant margin-free description of the dependence structure
of a random variable. By the famous theorem of \citet{sklar1959}, it
follows that the
distribution function $H$ of a bivariate random variable $(X,Y)$
can be represented in terms of the marginal distributions $F$ and $G$
of $X$ and~$Y$, that is,
\[
H (x,y) = C(F(x), G(y)),
\]
where $C$ denotes the copula, which characterizes the dependence between~$X$
and $Y$. Extreme-value copulas arise naturally as the possible
limits of
copulas of component-wise maxima of independent, identically
distributed or strongly mixing stationary sequences [see
\citet{deheuvels1984} and \citet{hsing1989}]. These copulas provide
flexible tools for modeling joint extremes in risk management.
An important application of extreme-value copulas appears in the
modeling of data with positive dependence, and in contrast to the more
popular class of Archimedean copulas they are not symmetric
[see \citet{tawn1988} 
or \citet{ghokhoriv1998}].
Further applications can be found in \citet{colheftaw1999} or
\citet{cebdenlam2003} among others. A copula $C$ is an extreme-value
copula if and only if it has a~representation of the form
%
%
\begin{equation}\label{11}
C(y^{1-t},y^t) = y^{A(t)}\qquad\forall y,t \in[0,1],
\end{equation}
where $A\dvtx[0,1] \to[1/2,1]$ is a convex function satisfying $\max\{
s,1-s\} \le A(s) \le1$, which is called
Pickands dependence function.
The representation of (\ref{11}) of the extreme-value copula $C$
depends only on the one-dimensional
function $A$ and statistical inference on a bivariate extreme-value
copula $C$ may now be reduced
to inference on its Pickands dependence function~$A$.

The problem of estimating Pickands dependence function nonparametrically
has found considerable attention in the literature.
Roughly speaking, there exist two classes of estimators.
The classical nonparametric
estimator is that of \citet{pickands1981} [see \citet
{deheuvels1991} for
its asymptotic properties] and several variants have been discussed.
Alternative estimators have been proposed and
investigated in the papers by \citet{capfougen1997}, \citet
{jimvilflo2001},
\citet{haltaj2000}, \citet{segers2007} and \citet
{zhawellpen2008},
where the last-named authors also discussed the
multivariate case. In most references, the estimators of Pickands dependence
function are constructed assuming knowledge of the marginal
distributions. Recently \citet{genest2009} proposed
rank-based versions of the estimators of \citet{pickands1981} and
\citet{capfougen1997}, which do not require
knowledge of the marginal distributions. In general, all of these
estimators are neither convex nor do they satisfy
the boundary restriction $\max\{t,1-t\} \le A(t) \le1$, in
particular the endpoint constrains $A(0)=A(1)=1$. However,
the estimators can be modified without changing their asymptotic
properties in such a way that these constraints are satisfied, see,
for example,~\citet{filsvilletard2008}.\looseness=-1

Before the specific model of an extreme-value copula is selected, it is
necessary to check this assumption
by a statistical test, that is a test for the hypotheses
%
%
\begin{equation}
\label{hypextr} H_0\dvtx C \in\mathcal{C} \quad\mbox{vs.}\quad H_1\dvtx
C\notin\mathcal{C},
\end{equation}
where $\mathcal{C}$ denotes the class of all copulas satisfying (\ref{11}).
Throughout this paper, we call (\ref{hypextr}) the hypothesis of
extreme-value dependence.\vadjust{\goodbreak}
The problem of testing this hypothesis has found much less attention in
the literature.
To our best knowledge, only two tests of extremeness are currently
available in the literature. The first one was proposed by \citet
{ghokhoriv1998}. It exploits the fact that for an extreme-value
copula the random variable $W=H(X,Y)=C(F(X),G(Y))$
satisfies the identity
%
%
\begin{equation} \label{HOW}
-1+8 \mathbb{E}[W]-9 \mathbb{E}[W^2 ]=0.
\end{equation}
The properties of this test have been studied by \citet{ghogennes2009},
who determined the finite- and large-sample variance of the test
statistic. In particular, the test proposed by \citet
{ghokhoriv1998} is not
consistent against alternatives satisfying (\ref{HOW}). The second
class of tests was recently introduced by \citet{koyan2010} who
proposed to
compare the empirical copula and a copula estimator which is
constructed from the
estimators proposed by \citet{genest2009} under the assumption of an
extreme-value copula. These tests are only consistent against
alternatives that are left tail decreasing in both arguments
and satisfy strong smoothness assumptions on the copula and convexity
assumptions on an analogue of
Pickands dependence function, which are hard to verify analytically.

The present paper has two purposes. The first is the development of
some alternative estimators of Pickands dependence function
using the principle of minimum distance estimation.
We propose
to consider the best approximation of the logarithm of the empirical
copula $\hat C$ evaluated in the point $(y^{1-t},y^t)$, that is, $\log
\hat C(y^{1-t},y^t)$, by functions of the form
%
%
\begin{equation}\label{12}
\log(y)A(t)
\end{equation}
with respect to a weighted $L^2$-distance. It turns out that the
minimal distance and the corresponding optimal function can be
determined explicitly.
On the basis of this result, and by choosing various weight functions
in the $L^2$-distance, we obtain an
infinite-dimensional class
of estimators for the function $A$.
Our approach is closely related to the theory of $Z$-estimation and
in Section \ref{secmindist} we indicate how this point of view
provides several interesting relationships
between the different concepts for constructing estimates of Pickands
dependence function.

The second purpose of the paper is to present a new test for the hypothesis
of extreme-value dependence, which is consistent against a much broader
class of alternatives than the
tests which have been proposed so far.
Here our approach is based on an estimator
of a weighted minimum $L^2$-distance between the true copula and the
class of
functions satisfying (\ref{12}) and the corresponding tests are
consistent with respect to all positive quadrant
dependent alternatives satisfying weak differentiability assumptions of
first order.
To our best knowledge, this method provides the first test in this context
which is consistent against such a general class of alternatives.
Moreover, in contrast to \citet{ghokhoriv1998} and \citet
{koyan2010} we
also provide a weak convergence result under
fixed alternative which can be used for studying the power of the test.

The remaining part of the paper is organized as follows. In Section
\ref{secevd},
we consider the approximation problem from a theoretical point of view.
In particular, we derive explicit representations for the minimal $L^2$-distance
between the logarithm of the copula and its best approximation by a
function of
the form (\ref{12}), which will be the basis for all statistical
applications in this
paper. The new estimators, say $\hat A_n$, are defined in Section \ref
{secmindist},
where we also prove weak convergence of the process $\{ \sqrt{n} (\hat
A_{n}(t)-A(t)) \}_{t \in[0,1]}$ in the space of uniformly bounded
functions on the
interval $[0,1]$ under appropriate assumptions on the weight function
used in the
$L^2$-distance. Furthermore, we give a theoretical and empirical
comparison of
the new estimators with the estimators proposed in \citet{genest2009}.
We will
also determine ``optimal'' estimators in the proposed class by
minimizing the asymptotic MSE with respect to the choice of the weight function
used in the $L^2$-distance. In particular, we demonstrate that some of
the new
estimators have a substantially smaller asymptotic variance than the estimators
proposed by the last-named authors. We also provide a simulation study
in order
to investigate the finite sample properties of the different estimates.
In Section \ref{sectestevd},
we introduce and investigate the new test of extreme-value dependence.
In particular,
we derive the asymptotic distribution of the test statistic under the
null hypothesis as
well as under the alternative. In order to approximate the critical
values of
the test, we introduce a multiplier bootstrap procedure, prove its
consistency and study
its finite sample properties by means of a simulation study. Finally,
most of the technical
details are deferred to the \hyperref[app]{Appendix}. 

\section{A measure of extreme-value dependence}\label{secevd}

Let $\mathcal{A}$ denote the set of all functions
$
A\dvtx[0,1]\rightarrow[1/2,1],
$
and define $\Pi$ as the copula corresponding to independent random
variables, that is, $\Pi(u,v) = uv$. Throughout this paper, we assume that
the copula $C$ satisfies $C\geq\Pi$ which holds for any extreme-value
copula due to the lower bound for the function $A$. As pointed out by
\citet{scaillet2005}, this property is equivalent to the concept
of positive
quadrant dependence, that is,
%
%
\begin{equation}
\mathbb{P}(X \leq x, Y \leq y) \geq \mathbb{P}(X \leq x) \mathbb{P}(Y \leq y)\qquad \forall (x,y)
\in\mathbb{R}^2.
\end{equation}
For a copula with this property, we define the weighted $L^2$-distance
%
%
\begin{equation} \label{21}
M_{h}(C,A) = \int_{(0,1)^2} \bigl( \log C(y^{1-t},y^t) - \log(y)A(t)
\bigr)^2 h(y) \,d(y,t),
\end{equation}
where $h\dvtx[0,1]\rightarrow\R^+$ is a continuous weight function.

The following result is essential for our approach and provides an explicit
expression for the best $L^2$-approximation of the logarithm of the copula
by the logarithm of a function of the form (\ref{11}) and as a by-product
characterizes the function $A^*$ minimizing $M_{h}(C,A) $.
\begin{theorem} \label{theobestapprox}
Assume that the given copula satisfies $C\geq\Pi^\kappa$ for some
$\kappa\geq1$
and that the weight function $h$ satisfies $\int_0^1 (\log y)^2 h(y)
\,dy < \infty$. Then
the function
\[
A^* =\arg\min\{M_{h}(C,A) | A \in\mathcal{A} \}
\]
is unique and given by
%
%
\begin{equation} \label{astar}
A^*(t) = B_h^{-1} \int_0^1 \frac{\log C(y^{1-t},y^t)}{\log y} h^*(y)
\,dy ,
\end{equation}
where the associated weight function $h^*$ is defined by
%
%
\begin{equation} \label{hstar}
h^*(y) = \log^2(y)h(y),\qquad y\in(0,1),
\end{equation}
and
%
%
\begin{equation} \label{22a}
B_h = \int_0^1 (\log y)^2 h(y) \,dy = \int^1_0 h^*(y) \,dy.
\end{equation}
Moreover, the minimal $L^2$-distance between the logarithms of the
given copula
and the class of functions of the form (\ref{12}) is given by
%
%
\begin{equation} \label{MCastar}
M_{h}(C,A^*) = \int_{(0,1)^2}\!\biggl( \frac{\log C(y^{1-t}, y^t)}{\log y}\biggr)^2 h^*(y) \,d(y,t)
- B_h \int^1_0\! (A^*(t))^2 \,dt.\hspace*{-32pt}
\end{equation}
\end{theorem}
\begin{pf}
Since\vspace*{1pt} $C\geq\Pi^\kappa$, we get $0 \geq\log C(y^{1-t},y^t) \geq
\kappa\log y$ and thus\break $|{\log C}(y^{1-t},y^t)| \leq\kappa|{\log y}|$
and all integrals exist.
Rewriting the $L^2$ distance in (\ref{21}) gives
\[
M_{h}(C,A) = \int_0^1 \int_0^{1} \biggl( \frac{\log
C(y^{1-t},y^t)}{\log y}- A(t) \biggr)^2 (\log y)^2 h(y) \,dy \,dt
\]
and the assertion is now obvious.
\end{pf}

Note that $A^*(t)=A(t)$ if $C$ is an extreme-value copula of the form
(\ref{11})
with Pickands dependence function $A$. Furthermore, the following lemma shows
that the minimizing function $A^*$ defined in (\ref{astar}) satisfies
the boundary
conditions of Pickands dependence functions.
\begin{lemma} \label{conditionpick}
Assume that $C$ is a copula
satisfying $C \ge\Pi$. Then the function $A^*$ defined
in (\ref{astar}) has the following properties:
\begin{longlist}[(iii)]
\item[(i)] $A^*(0)=A^*(1)=1$,
\item[(ii)] $A^*(t)\geq t \vee(1-t)$,
\item[(iii)] $A^*(t)\leq1$.
\end{longlist}
\end{lemma}
\begin{pf}
Assertion (i) is obvious. For a proof of (ii),
one uses the Fr\'{e}chet--Hoeffding bound $C(u,v)\leq u\wedge v$ [see,
e.g., \citet{nelsen2006}]
and obtains the assertion by a direct calculation.
Similarly, assertion
(iii) follows from the inequality $C\geq\Pi$.
\end{pf}

Unfortunately, the function $A^*$ is in general not convex for every
copula satisfying $C\geq\Pi$.
A counterexample can be derived from Theorem 3.2.2 in \citet
{nelsen2006} and is given by the following shuffle of
the copula $u\wedge v$:
%
%
\begin{equation}\label{ceconvex}
C(u,v)= \cases{
\min\{u,v,1/2\}, \qquad\mbox{$(u,v)\in\bigl[0,\sqrt{1/2}\bigr]^2$}, \vspace*{2pt}\cr
\min\bigl\{u,v+1/2-\sqrt{1/2}\bigr\}, \vspace*{2pt}\cr
\hspace*{91pt}(u,v)\in\bigl[0,\sqrt{1/2}\bigr]\times\bigl[\sqrt{1/2},1\bigr],\vspace*{2pt}\cr
\min\bigl\{u+1/2-\sqrt{1/2},v\bigr\}, \vspace*{2pt}\cr
\hspace*{91pt}(u,v)\in \bigl[\sqrt{1/2},1\bigr]\times\bigl[0,\sqrt{1/2}\bigr],\vspace*{2pt}\cr
\min\bigl\{u,v,u+v+1/2-2\sqrt{1/2}\bigr\}, \vspace*{2pt}\cr
\hspace*{91pt}\mbox{$(u,v)\in\bigl[\sqrt{1/2},1\bigr]^2$},}
\end{equation}
for which an easy calculation shows that the mapping $t\,{\mapsto}\,-\log C(1/2^{1-t},\!1/2^t)$ is not convex.
Consequently, one can find a weight function $h$ such that the
corresponding best approximating function
$A^*$ is not convex.

With the notation
%
%
\begin{equation} \label{fy}
f_y(t)=C(y^{1-t},y^t),
\end{equation}
the function $A^*$ is convex (for every weight function $h$) if and
only if the function
$g_y(t)= -\log f_y(t)$ is convex for every $y\in(0,1)$.
The following lemma is now obvious.
\begin{lemma} \label{convex}
If the function $t \rightarrow f_y(t) =C(y^{1-t},y^t) $ is twice
differentiable and the inequality
\[
[f_y'(t)]^2 \geq f_y''(t) f_y(t)
\]
holds for every $(y,t)\in(0,1)^2$, then the best approximation $A^*$
defined by~(\ref{astar}) is convex.
\end{lemma}

It is worthwhile to mention that the function $A^*$ is convex for some
frequently considered classes of copulas,
which will be illustrated in the following examples.
\begin{example} \label{claytonastar}
Consider the Clayton copula
%
%
\begin{equation} \label{copclay}
C_{\mathrm{Clayton}}(u,v;\theta) = ( u^{-\theta} + v^{-\theta
} - 1)^{-1/\theta},\qquad \theta>0.
\end{equation}
Then a tedious calculation yields
\begin{eqnarray*}
&& [f_y'(t)]^2-f_y''(t)f_y(t) \\
&&\qquad= \theta\log^2(y) \{ C_{\mathrm
{Clayton}}(y^{1-t},y^t;\theta)\} ^{2+2\theta}\bigl(
4y^{-\theta}-y^{-\theta t}-y^{-\theta(1-t)} \bigr)\\
&&\qquad\geq \theta\log^2(y) \{ C_{\mathrm
{Clayton}}(y^{1-t},y^t;\theta)\}^{2+2\theta} (
3y^{-\theta}- 1) \geq0,
\end{eqnarray*}
where the inequalities follow observing that $m(t)=y^{-\theta
t}+y^{-\theta(1-t)} \leq m(0)=1+y^{-\theta}$ and $y^{-\theta}\geq1$.
Therefore, we obtain from Lemma \ref{convex} that the best
approximation $A^*$ is convex and corresponds to an
extreme-value copula.
\end{example}
\begin{example} \label{exhalphahk}
In the following, we discuss the weight function $h_k(y)= -y^k / \log
y$ ($k\geq0$) with
associated function $h_k^*(y)= -y^k\log y $,
which will be used later for the construction of the new estimators of
Pickands dependence function.
On the one hand this
choice is made for mathematical convenience, because it allows an
explicit calculations of the asymptotic variance $A^*$
in specific examples. On the other hand, estimates constructed on the
basis of this weight function turn out to have good
asymptotic and finite sample properties (see the discussion in Section
\ref{simulation}).
It follows that
\[
B_{h_k}= - \int_0^1 y^k \log y \,dy=(k+1)^{-2}
\]
and
%
%
\begin{equation} \label{27}
A^*(t)=-(k+1)^2 \int_0^1 \log C(y^{1-t},y^t) y^k \,dy,
\end{equation}
which simplifies in the case $k=0$ to the representation
%
%
\begin{equation} \label{26}
A^*(t)= - \int^1_0 \log C(y^{1-t},y^t) \,dy.
\end{equation}
\end{example}
\begin{example}
$\!\!\!$In the following, we calculate the minimal distance $M_{h}(C,\allowbreak A^*)$
and its corresponding best approximation $A^*$ for two copula families
and the associated weight function $h_1^*(y) = -y\log y$
from Example \ref{exhalphahk}. First, we investigate the Gaussian
copula defined by
\[
C_\rho(u,v)=\Phi_2(\Phi^{-1}(u), \Phi^{-1}(v), \rho),
\]
where $\Phi$ is the standard normal distribution function and $\Phi
_2(\cdot,\cdot,\rho)$ is the distribution function of a bivariate
normal random variable with standard normally distributed margins and
correlation $\rho\in[0,1]$. For the limiting cases $\rho=0$ and
$\rho=1$, we obtain the independence and perfect dependence copula,
respectively, while for $\rho\in(0,1)$ the copula $C_\rho$ is not an
extreme-value copula.
The minimal distances are plotted as a function of $\rho$ in the left
part of the first line of Figure \ref{picdist}.
In the right part, we show some functions $A^*$ corresponding to the
best approximation of the logarithm of the Gaussian copula by
a function of the form (\ref{12}).
We note that all functions $A^*$ are convex although $C_\rho$ is only
an extreme value copula in the case $\rho=0$.

%
\begin{figure}

\includegraphics{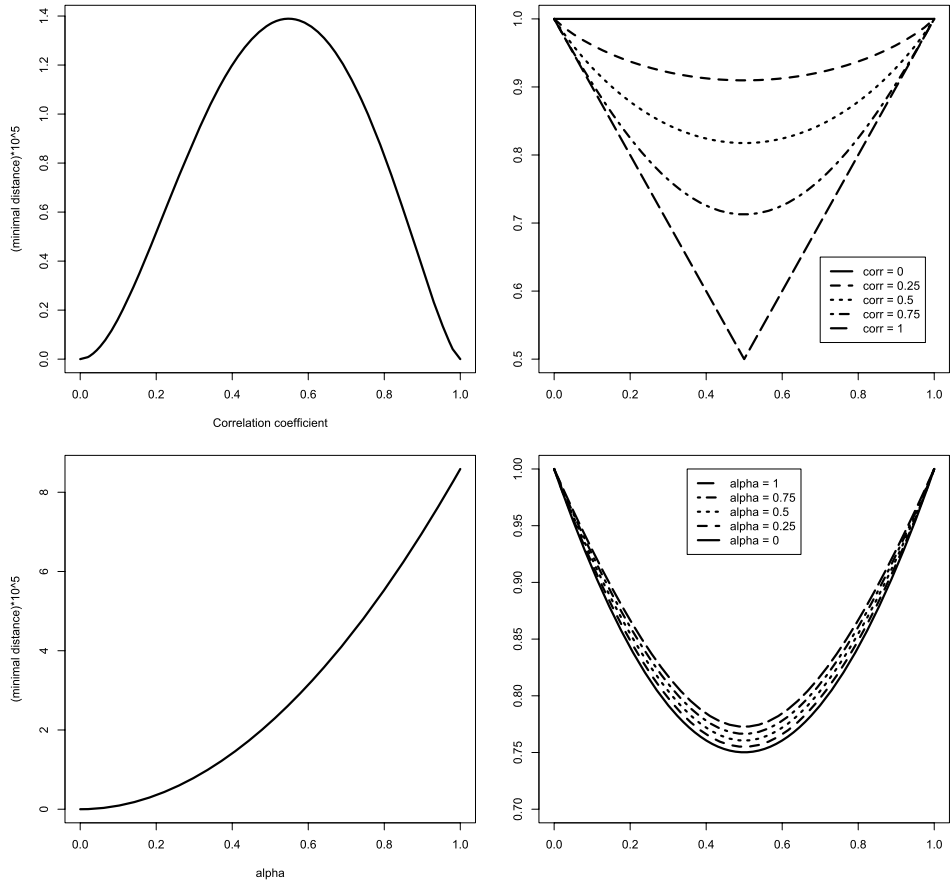}

\caption{Left: minimal distances $M_{h}(C,A^*)\times10^5$
for the Gaussian copula (as a function of its correlation coefficient)
and for the convex combination of a Gumbel and a Clayton copula (as a
function of the parameter
$\alpha$ in the convex combination). Right: the functions~$A^*$ corresponding
to the best approximations by functions of the form (\protect\ref
{12}).} \label{picdist}
\vspace*{3pt}
\end{figure}

In the second example, we consider a convex combination of a Gumbel
copula with parameter
$\theta_1=\log2/\log1.5$ (corresponding to a coefficient of tail
dependence of $0.5$) and a Clayton copula with parameter $\theta_2=2$,
that is,
\[
C_\alpha(u, v)= \alpha C_{\mathrm{Clayton}}(u,v;\theta
_2)+(1-\alpha) C_{\mathrm{Gumbel}}(u,v;\theta_1),\qquad\alpha\in
[0,1],
\]
where the Clayton copula is given in (\ref{copclay}) and the
Gumbel copula is defined by
\[
C_{\mathrm{Gumbel}}(u,v;\theta) = \exp\bigl( - \{ (-\log
u)^\theta+ (-\log v)^\theta \}^{1/\theta} \bigr),\qquad
\theta>1.
\]
Note that only the Gumbel copula is an extreme-value copula and
obtained for \mbox{$\alpha=0$}.
The minimal distances are depicted in the left part of the lower panel
of Figure~\ref{picdist}
as a function of $\alpha$. In the right part, we show the functions~$A^*$ corresponding to the best approximation
of the logarithm of $C_\alpha$ by a~function of the form (\ref{12}).
Again all approximations are convex, which means that $A^*$ corresponds
in fact to an extreme value copula.
\end{example}

\section{A class of minimum distance estimators}\label{secmindist}

\subsection{Pickands and CFG estimators}

Let $(X_1,Y_1),\ldots,(X_n,Y_n)$ denote\break a~sample of independent identically
distributed bivariate random variables with copula $C$ and marginals
$F$ and
$G$. Most of the estimates which have been proposed in the literature
so far are based on the
fact that the random variable
\[
\xi(t)= \frac{- \log F(X)}{1-t} \wedge\frac{- \log G(Y)}{t}
\]
is exponentially distributed with parameter $A(t)$. In particular, we
have $E[\xi(t) ] = 1/ A(t)$.
If the marginal distributions would be known, an estimate of $A(t)$
could be obtained
by the method of moments. In the case of unknown marginals, \citet
{genest2009}
proposed to replace $F$ and $G$ by their empirical counterparts and obtained
\[
\hat A^P_{n,r}(t)= \Biggl( \frac{1}{n} \sum^n_{i=1} \hat\xi_i (t)
\Biggr)^{-1}
\]
as a rank-based version of Pickands estimate, where
\[
\hat\xi_i (t)= \frac{- \log\hat F_n (X_i)}{1-t} \wedge\frac{-
\log\hat G_n(Y_i) }{t},\qquad i=1,\ldots,n,
\]
and
%
%
\begin{equation} \label{31a}
\qquad\hat F_n (X_i)\,{=}\,\frac{1}{n+1} \sum^n_{j=1} I \{X_j\,{\leq}\,X_i \} \quad\mbox{and}\quad
\hat G_n(Y_i)\,{=}\,\frac{1}{n+1} \sum^n_{j=1} I \{Y_j\,{\leq}\,Y_i \}
\end{equation}
denote the (slightly modified) empirical distribution functions of the
samples $\{X_j\}^n_{j=1}$ and $\{Y_j\}^n_{j=1}$ at the points $X_i$ and
$Y_i$, respectively.
Similarly, observing the identity
$E[\log\xi(t) ] = - \log A(t) - \gamma$ (here $\gamma=- \int
^\infty_0 \log x e^{-x}\,dx$ denotes Euler's
constant),\vadjust{\goodbreak}
they obtained a rank-based version of the estimate proposed by
\citet
{capfougen1997}, that is,
\[
\hat A^{\mathrm{CFG}}_{n,r} (t)= \exp\Biggl( - \gamma- \frac{1}{n} \sum
^n_{i=1} \log\hat\xi_i (t) \Biggr).
\]
For illustrative purposes, we finally recall two integral
representations for the rank-based version of
Pickands and CFG estimate,
which we use in~Sec\-tion~\ref{subsecmundzest} to put all estimates
considered in this paper in a general
context, that is,
%
%
\begin{eqnarray}
\label{pickdef}
\frac{1}{\hat A^P_{n,r}(t)} &=& \int^1_0 \frac{\hat C_n (y^{1-t},
y^t)}{y}\,dy, \\[-3pt]
\label{cfgdef}
\gamma+ \log\hat A^{\mathrm{CFG}}_{n,r} (t) &=& \int^1_0 \frac{\hat C_n
(y^{1-t}, y^t) - I \{ y > e^{-1} \}}{\log y}\,dy,
\end{eqnarray}
where
%
%
\begin{equation}
\label{31}
\hat C_n(u,v)= \frac{1}{n} \sum^n_{i=1} I \{ \hat F_n (X_i) \leq u,
\hat G_n(Y_i) \leq v \}
\end{equation}
denotes the empirical copula and $ \hat F_n (X_i)$, $\hat G_n(Y_i) $
are defined in (\ref{31a})
[see \citet{genest2009} for more details].

\vspace*{-3pt}\subsection{New estimators and weak convergence} \label{subsecnewest}

Theorem \ref{theobestapprox} suggests to define a class of new
estimators for
Pickands dependence function by replacing the unknown copula in (\ref{astar})
through the empirical copula defined in (\ref{31}). The asymptotic
properties of
the corresponding estimators will be investigated in this section. For technical
reasons, we require that the argument in the logarithm in the representation
(\ref{astar}) is positive and propose to use the estimator
%
%
\begin{equation} \label{tildeCn}
\tilde C_n = \hat C_n \vee n^{-\gamma},
\end{equation}
where the constant $\gamma$ satisfies $\gamma>1/2$ and the empirical copula
$\hat C_n$ is defined in~(\ref{31}).

For the subsequent proofs, we will need a result on the weak
convergence of the empirical copula process with estimated margins.
While this problem has been considered by many authors
[see, e.g., \citet{ruschendorf1976}, \citet{ferradweg2004} or
\citet
{tsukahara2005} among others], all of them assume that the copula
has continuous partial derivatives on the whole unit square $[0,1]^2$.
However, as was pointed out by \citet{segers2010}, there is
only one extreme-value copula that has this property. Luckily,
in a~remarkable paper \citet{segers2010} was able to show that the following
condition is
sufficient for weak convergence of the empirical copula process
%
%
\begin{equation} \label{asproc}
\partial_j C\mbox{ exists and is continuous on } \{(u_1,u_2)\in
[0,1]^2 | u_j \in(0,1)\}
\end{equation}
($j=1,2$). This condition can be shown to hold for any extreme-value
copula with continuously differentiable Pickands function $A$ [see
\citet{segers2010}].\vadjust{\goodbreak} Moreover, under this assumption, the process
$\sqrt{n}(\tilde C_n - C)$ shows the same
limiting behavior as the empirical copula process $\sqrt{n}({\hat C}_n
- C)$, that is,
%
%
\begin{equation}
\label{32}
\sqrt{n}(\tilde C_n -C) \weak\G_C,
\end{equation}
where the symbol $\weak$ denotes weak convergence in $l^\infty[0,1]^2$.
%
Here, $\G_C$ is a~Gaussian field on the square $[0,1]^2$ which admits
the representation
\[
\G_C(\x) = \B_C(\x) - \partial_1 C(\x)\B_C(x_1,1) - \partial
_2C(\x) \B_C(1,x_2),
\]
where $\x= (x_1,x_2), \B_C$ is a bivariate pinned $C$-Brownian sheet
on the square
$[0,1]^2 $ with covariance kernel given by
\[
\Cov(\B_C(\x),\B_C(\y)) = C(\x\wedge\y) - C(\x) C(\y),
\]
and the minimum $\x\wedge\y$ is understood component-wise. Observing
the representation (\ref{astar}), we obtain the estimator
%
%
\begin{equation} \label{33}
\hat A_{n,h}(t) = B_h^{-1} \int_0^1 \frac{ \log\tilde
C_n(y^{1-t},y^t)}{\log y} h^*(y) \,dy
\end{equation}
for Pickands dependence function, where $\tilde C_n$ is defined in
(\ref{tildeCn}).
Note that this relation specifies an infinite-dimensional class of
estimators indexed
by the set of all admissible weight functions. The following results specify
the asymptotic properties of these estimators. We begin with a slightly
more general
statement, which shows weak convergence for the weighted integrated process\vspace*{-4pt}
\[
\sqrt{n}\W_{n,w}(t) = \sqrt n \int_0^1 \log\frac{\tilde
C_n(y^{1-t}, y^t)}{C(y^{1-t}, y^t)} w(y,t) \,dy,
\]
where the weight function
$ w\dvtx[0,1]^2 \rightarrow\bar\R$ depends on $y$ and $t$. The result
(and some arguments in its proof) are also needed in Section
\ref{sectestevd}.\vspace*{-6pt}
\begin{theorem}\label{theoh1}
Assume that for the weight function $w\dvtx[0,1]^2 \rightarrow\bar\R$
there exists a function $\bar w\dvtx[0,1] \rightarrow\bar\R_0^+$ such that
%
%
\begin{eqnarray}
\label{cw1}
\forall (y,t) \in[0,1]^2\qquad|w(y,t)|&\leq&\bar w(y),\\[-2pt]
\label{cw2}
\forall\varepsilon>0\qquad\sup_{y \in[\varepsilon,1]} \bar w(y) &<&
\infty,\\[-2pt]
\label{cw3}
\int_0^1 \bar w(y)y^{-\lambda}\,dy&<&\infty
\end{eqnarray}
for some $\lambda>1 $. {If the copula $C$ satisfies (\ref{asproc})}
and $C \geq\Pi$, then
we have for any $\gamma\in(1/2, \lambda/2)$ as $n\rightarrow\infty$
%
%
\begin{eqnarray} \label{wk1}
&&\sqrt{n}\W_{n,w}(t) = \sqrt n \int_0^1 \log\frac{\tilde
C_n(y^{1-t}, y^t)}{C(y^{1-t}, y^t)} w(y,t) \,dy \nonumber\\[-10pt]\\[-10pt]
&&\hspace*{4pt}\quad\weak\quad\W_{C,w}(t) = \int_0^1 \frac{\G_C(y^{1-t},
y^t)}{C(y^{1-t},y^t)} w(y,t) \,dy \nonumber
\end{eqnarray}
in $ l^\infty[0,1]$.\vadjust{\eject}
\end{theorem}

The following result is now an immediate consequence of Theorem \ref{theoh1}
using $w(y,t) := -B_h^{-1}h^*(y)$ [recall the definition of the associated
weight function $h^*$ in (\ref{hstar})]
and yields the weak
convergence of the process \mbox{$\sqrt{n}(\hat A_{n,h} - A^*)$} for a broad
class of weight functions.
\begin{theorem}\label{theopickest1}
If the copula $C \geq\Pi$ satisfies condition (\ref{asproc}) and the
weight function $h$ satisfies
the conditions
%
%
\begin{eqnarray}
\label{20a} \mbox{for all } \eps>0\qquad\sup_{y\in[\eps,1]} \biggl|
\frac{h^*(y)}{\log y} \biggr| &<&\infty,\\
\label{20c} \int_0^1 h^*(y) (-\log y)^{-1} y^{-\lambda} \,dy&<&\infty
\end{eqnarray}
for some $\lambda>1$,
then we have for any $\gamma\in(1/2,\lambda/2)$ as $n\rightarrow
\infty$
\[
\A_{n,h} = \sqrt{n}(\hat A_{n,h} - A^*) \quad\weak\quad\A_{C,h} \mbox{ in }
l^\infty[0,1],
\]
where the process $\A_{C,h}$ is given by
%
%
\begin{equation} \label{34}
\A_{C,h}(t) = B_h^{-1} \int_0^1 \frac{\G_C(y^{1-t},
y^t)}{C(y^{1-t},y^t)} \frac{h^*(y)}{\log y} \,dy.
\end{equation}
\end{theorem}
\begin{remark}
(a) Conditions (\ref{20a}) and (\ref{20c}) restrict the behavior of
the function $h^*$ near the boundary
of the interval $[0,1]$.
A simple sufficient condition for (\ref{20a}) and (\ref{20c}) is
given by
\[
\sup_{x\in[0,1]} \biggl| \frac{h^*(x)}{x^\alpha(1-x)^\beta}
\biggr| < \infty
\]
for some $\alpha>0,\beta\geq1$. In this case, $\lambda$ can be
chosen as $1+\alpha/2$.

(b) In the construction discussed so far, it is also possible to use
weight functions that depend on
$t$, that is, functions of the form $\tilde h^*(y,t)$. As long as
$\tilde
h^*(y,t) >0$ for $(y,t) \in(0,1)\times[0,1]$,
the corresponding best approximation $A^*$
will still be well defined and correspond to the Pickands dependence
function if $C$ is an extreme-value copula.
Theorem \ref{theoh1} provides the asymptotic properties of the
corresponding estimator $A$ if we set $w(y,t) := \tilde h^*(y,t)/(-\log y)$
and assume that $\int_0^1 \tilde h^*(y,t)\,dy = 1$ for all $t$. However,
for the sake of a clear presentation, we will only use weight functions
that do not depend on $t$.
\end{remark}

Note that Theorem \ref{theopickest1} is also correct if the given
copula is not
an extreme-value copula. In other words: it establishes weak
convergence of
the process $\sqrt{n}(\hat A_{n,h}- A^*)$ to a centered Gaussian process,
where $A^*$ denotes the function corresponding to the best
approximation of
the logarithm of the copula $C$ by a function of the form (\ref{12}).
If $A^*$
is convex, it corresponds to an extreme-value copula and coincides with Pickands
dependence function. Note also that Theorem~\ref{theopickest1}
excludes the
case $h_0^*(y) = - \log y$, because condition (\ref{20c}) is not
satisfied for this
weight function. Nevertheless, under the additional assumption that $C$
is an
extreme-value copula with twice continuously differentiable Pickands dependence
function $A$, the assertion of the preceding theorem is still valid.
\begin{theorem}\label{theopickest2}
Assume that
$C$ is an extreme-value copula with twice continuously differentiable
Pickands dependence function $A$. For the weight function
$h^*_0(y)=-\log y$, we have for any $\gamma\in(1/2,3/4)$ as
$n\rightarrow\infty$
\[
\A_{n,h_0}(t) = \sqrt{n}(\hat A_{n,h_0 } - A)(t) \quad\weak\quad\A_{C,h_0
}(t) = - \int_0^1 \frac{\G_C(y^{1-t}, y^t)}{C(y^{1-t},y^t)} \,dy
\]
in $ l^\infty[0,1]$, where $\hat A_{n,h_0 }(t) =- \int_0^1 \log
\tilde C_n(y^{1-t},y^t) \,dy$.
\end{theorem}
\begin{remark}
(a) If the marginals of $(X,Y)$ are independent the distribution of
the random variable $\A_{\Pi,h_0}$ coincides with
the distribution of the random variable $\A^{P}_r=- \int_0^1 {\G_\Pi
(y^{1-t}, y^t)}{y}^{-1} \,dy$, which appears as the weak limit of the
appropriately standardized Pickands estimator; see \citet{genest2009}.
In fact, a much more general
statement is true: by using weight functions $\tilde h^*(y,t)$
depending on $t$
it is possible to obtain
for any extreme-value copula estimators of the form (\ref{33})
which show the same limiting behavior as the estimators proposed by
\citet{genest2009}.
This already indicates that for any extreme-value copula it is possible
to find weight
functions which will make the new minimum distance estimators
asymptotically at least as
efficient (in fact better, as will be shown in Section \ref
{subsecoptweight}) as the estimators
introduced by \citet{genest2009}.

(b) A careful inspection of the proof of Theorem \ref{theoh1} reveals
that the condition $C\geq\Pi$ can be relaxed to $C\geq\Pi^\kappa$
for some $\kappa>1$, if one imposes stronger conditions on the weight function.

(c) The estimator depends on the parameter $\gamma$ which is used for
the construction of the
statistic $\tilde C_n = \hat C_n \vee n^{-\gamma} $. This modification
is only made for technical purposes and
from a practical point of view the behavior of the estimators does
not change substantially provided that $\gamma$ is chosen larger than
$ 2/3$.
\end{remark}
\begin{remark} \label{bayes}
The new estimators can be alternatively motivated observing that the
identity (\ref{11}) yields the representation $A(t)\,{=}\,\log
C(y^{1-t}, y^t)/\allowbreak \log y$ for any $y \in(0,1)$. This leads to a simple
class of estimators, that is,
\[
\tilde A_{n,\delta_y} (t)= \frac{\log\tilde C_n(y^{1-t},y^t)}{\log
y};\qquad y \in(0,1),
\]
where $\delta_y$ is the Dirac measure at the point $y$ and $\tilde
C_n$ is defined in (\ref{tildeCn}). By averaging these estimators with
respect to a distribution,
say $\pi$, we obtain estimators of the form
\[
\tilde A_{n,\pi} (t)= \int^1_0 \frac{\log\tilde
C_n(y^{1-t},y^t)}{\log y} \pi(dy),
\]
which coincide with the estimators obtained by the concept of best
$L^2$-ap\-proximation.\vspace*{-3pt}
\end{remark}

\subsection{A special class of weight functions} \label{exampleexp}

In this subsection, we illustrate the results investigating Example
\ref{exhalphahk} discussed at the end of
Section \ref{secevd}. For the associated weight function
$h_k^*(x)=-y^k\log y$ with $k\geq0$, we obtain
%
%
\begin{equation}
\label{36}
\hat A_{n,h_k}(t) =-(k+1)^2 \int_0^1 \log\tilde C_n(y^{1-t},y^t)
y^k \,dy.
\end{equation}
The process $\{\A_{n,h_k}(t) \}_{t\in[0,1]}$ converge weakly in $
l^\infty[0,1]$ to the
process\break $\{ \A_{C,h_k}\}_{t\in[0,1]}$, which is given by
%
%
\begin{equation}
\A_{C,h_k}(t) =
-(k+1)^2 \int_0^1 \frac{\G_C(y^{1-t},y^t)}{C(y^{1-t},y^t)} y^k \,dy.
\end{equation}
Consequently, for $C\in\mathcal{C}$, the asymptotic variance
of $\hat A_{n,h_k}$ is obtained as
%
%
\begin{equation}\label{varhk}
\Var(\A_{C,h_k} (t) ) = (k+1)^4 \int_0^1\int_0^1 \sigma(u,v;t)
(uv)^{k-A(t)} \,du \,dv,
\end{equation}
where the function $\sigma$ is given by
\[
\sigma(u,v;t)=\Cov(\G_C(u^{1-t},u^t),\G_C(v^{1-t},v^t)).
\]
In order to find an explicit expression for these variances, we assume
that the function $A$ is differentiable and introduce the notation
\[
\mu(t)=A(t)-tA'(t),\qquad\nu(t)=A(t)+(1-t)A'(t),
\]
where $A'$ denotes the derivative of $A$.
The following results can be shown by similar arguments as given in
\citet{genest2009}; for
details, see \citet{techrep}.\vspace*{-3pt}

\begin{proposition} \label{propvarAk}
For $t\in[0,1]$, let $\bar\mu(t)=1-\mu(t)$ and $\bar\nu(t)=1-\nu
(t)$. If $C$ is an extreme-value copula
with Pickands dependence function $A$, then
the variance of the random variable $\A_{C,h_k} (t) $ is given by
\begin{eqnarray*}
\hspace*{-4pt}&& (k+1)^2\biggl\{ \frac{2(k+1)}{2k+2-A(t)} - \bigl(\mu(t)+\nu(t)-1\bigr)^2 \\
\hspace*{-4pt}&&\hphantom{(k+1)^2\biggl\{}{} - \frac{2\mu(t)\bar\mu(t) (k+1)}{2k+1+t} - \frac{2\nu
(t)\bar\nu(t) (k+1)}{2k+2-t}\\
\hspace*{-4pt}&&\hphantom{(k+1)^2\biggl\{}{} + 2\mu(t)\nu(t) \frac{(k+1)^2}{(1-t)t} \int_0^1 \biggl(
A(s)+(k+1)\biggl( \frac{1-s}{1-t} + \frac{s}{t}\biggr) - 1
\biggr)^{-2} \,ds \\
\hspace*{-4pt}&&\hphantom{(k+1)^2\biggl\{}{} - 2 \mu(t)\frac{(k+1)^2}{(1-t)t} \int_0^t \biggl( A(s) +
(k+t)\frac{1-s}{1-t} + \bigl(k+1-A(t)\bigr) \frac{s}{t} \biggr)^{-2} \,ds\\
\hspace*{-4pt}&&\hphantom{(k+1)^2\biggl\{}{} - 2 \nu(t)\frac{(k+1)^2}{(1-t)t} \int_t^1 \biggl( A(s) +
\bigl(k+1-A(t)\bigr)\frac{1-s}{1-t} \\
\hspace*{-4pt}&&\qquad\quad\hspace*{56.4pt}\hspace*{156pt}{} +(k+1-t)\frac{s}{t} \biggr)^{-2} \,ds
\biggr\}.
\end{eqnarray*}
\end{proposition}

Note that the limiting process in (\ref{34}) is a centered Gaussian process.
This means that, asymptotically, the quality of the new estimators [as
well as of the
estimators of \citet{genest2009}, which show a similar limiting
behavior] is
determined by the variance. Based on these observations, we will now
provide an asymptotic comparison of the new estimators
$ \hat A_{n, h_k}(t)$ with the estimators investigated
by \citet{genest2009}.
Some finite sample results will be presented in the following section
for various families of copulas.
For the sake of brevity, we restrict ourselves to the independence
copula~$\Pi$,
for which $A(t)\equiv1$.
In the case $k=0$, we obtain from Proposition~\ref{propvarAk} the same
variance as for the
rank-based version of Pickands estimator, that is,
\[
\Var(\A_{\Pi,h_0}) = \frac{3t(1-t)}{(2-t)(1+t)} = \Var( \A^P_r)
\]
[see Corollary 3.4 in \citet{genest2009}]
while the case $k > 0$ yields
\[
\Var(\A_{\Pi,h_k}) = \frac{(3+4k)(k+1)^2}{2k+1} \frac
{t(1-t)}{(2k+2-t)(2k+1+t)}.
\]
Investigating the derivative in $k$, it is easy to see that $\Var(\A
_{\Pi,h_k})$ is strictly decreasing in $k$ with
\[
\lim_{k\rightarrow\infty} \Var(\A_{\Pi,h_k}) = \frac{t(1-t)}{2}.
\]
Therefore, we have
\[
\Var( \A^P_r)=\Var(\A_{\Pi,h_0}) \geq\Var(\A_{\Pi,h_k})
\]
for all $k \geq0$ with strict inequality for all $k > 0$. This means
that for the independence copula all estimators
obtained by our approach with associated weight function
$h_k^*(y)=-y^k\log y$, $k >0$, have a smaller asymptotic variance than
the rank-based version of Pickands estimator. On the other hand, a
comparison with the
CFG estimator proposed by \citet{genest2009}
does not provide a clear picture about the superiority of one estimator
and we defer this comparison to the following
section, where optimal weight functions for the new estimates $ \hat
A_{n,h}$ are introduced.

\subsection{Optimal weight functions} \label{subsecoptweight}

In this section, we discuss asymptotically optimal weight functions
corresponding to the class of
estimates introduced in Section~\ref{subsecnewest}. As pointed out in
the previous section,
from an asymptotic point of view the mean squared error of the
estimates is dominated by the
variance and therefore we concentrate on weight functions minimizing
the asymptotic variance
of the estimate $ \hat A_{n,h}$. The finite sample properties of the
mean squared error of the
various estimates will be investigated by means of a simulation study
in Section \ref{simulation}.

Note that an optimal weight function depends on the point $t$ where
Pickands dependence function has to be estimated and on the unknown
copula. Therefore, an estimator with
an optimal weight function cannot be implemented in
concrete applications without preliminary knowledge about the copula.
However, it can serve as a benchmark for user-specified weight
functions. To be precise, observe that by Theorem
\ref{theopickest1} the variance of the limiting process $ \A_{C,h}$
is of the form
%
%
\begin{equation} \label{asymvarall}
V( \xi) = \int_0^1\int_0^1 k_t(x,y) \,d \xi(x) \,d \xi(y),
\end{equation}
where $\xi$ denotes a probability measure on the interval $[0,1]$
defined by $d \xi(x) = B_h^{-1} h^*(x)\,dx$
and the kernel $ k_t(x,y) $ is given by
\[
k_t(x,y) =
E \biggl[ \frac{\G_C(x^{1-t}, x^t)}{C(x^{1-t},x^t)\log x} \frac{\G
_C(y^{1-t}, y^t)}{C(y^{1-t},y^t)\log y} \biggr].
\]
It is easy to see that $V$ defines a convex function on the space of
all probability measures
on the interval $[0,1]$ and the existence of a minimizing measure follows
if the kernel $k_t$ is continuous on $[0,1]^2$. The following result
characterizes the minimizer of $V$
and is proved in the \hyperref[app]{Appendix}.
\begin{theorem}\label{effic}
A probability measure $\eta$ on the interval
$[0,1]$ minimizes~$V$ if and only if the inequality
%
%
\begin{equation} \label{lowbound}
\int_0^1 k_t(x,y) \,d \eta(y) \ge
\int_0^1\int_0^1 k_t(x,y) \,d \eta(x) \,d \eta(y)
\end{equation}
is satisfied for all ${x\in[0,1]} $.
\end{theorem}

Theorem \ref{effic} can be used to check the optimality of a given
weight function. For example,
if the copula $C$ is given by the independence copula $\Pi$ we have
\[
k_t(x,y) = \frac{(x^t \wedge y^t - (xy)^t) (x^{1-t} \wedge y^{1-t} -
(xy)^{1-t})}{x \log x\, y \log y},
\]
and it is easy to see that none of the associated weight functions
$h_k^*(y)=-y^k\log y$ with $k\geq0$
is optimal in the sense that it minimizes the asymptotic variance of
the estimate $\hat A_{n,h}$
with respect to the choice of the weight function.
On the other hand, the result is less useful for an explicit
computation of optimal weight functions.
Deriving an analytical expression for the optimal weight function seems
to be impossible,
even for the simple case of the independence copula.

However, approximations to the optimal weight function can
easily be computed numerically.
To be precise we approximate the double integral appearing in the
representation of $\Var(\A_{C,h}(t))$ by the finite sum
%
%
\begin{equation} \label{diskret}
V (\xi) \approx\sum_{i=1}^N \sum_{j=1}^N \xi_{i,N} \xi_{j,N}
k_t(i/N,j/N) = \Xi^T K_t \Xi,
\end{equation}
where\vspace*{1pt} $N\!\in\!\Nat$, $K_t\,{=}\,( k_t(i/N,j/N))_{i,j=1}^N$ denotes an
$N\!\times\!N$ matrix, $\Xi\,{=}\,(\xi_{i,N})_{i=1}^N$ is an vector of length $N$
and $\xi_{i,N} = \xi((i-1)/N,i/N]$ represents the mass of $\xi$
allocated to the interval $((i-1)/N,i/N]$ ($i=1,\ldots,N$).
Minimizing the right-hand side of the above equation with respect to
$\Xi$ under the constrains $ \xi_{i,N} \geq0$, $\sum_{i=1}^N \xi
_{i,N} = 1$
is a quadratic (convex) optimization problem which can be solved by
standard methods; see, for example,
\citet{nocedal2006} and approximations of the optimal weight function
can be calculated with arbitrary precision by increasing $N$.

\begin{figure}

\includegraphics{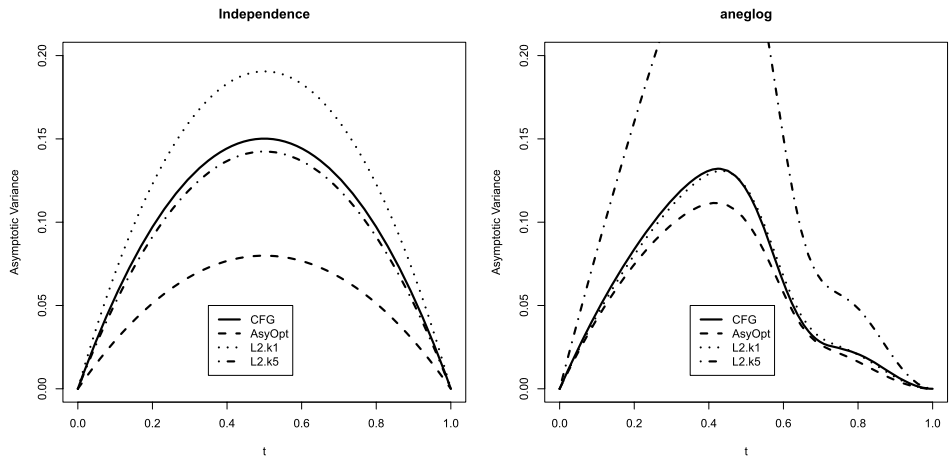}

\caption{Asymptotic variances of various estimators of the Pickands
dependence function. Left panel:
independence copula; right panel: asymmetric negative logistic model.}
\label{asymvar}
\end{figure}

In the remaining part of this section, we will compare the asymptotic variance
of the Pickands-, the CFG-estimator proposed by \citet{genest2009} and
the new estimates, where the new estimators are based on the weight
functions $h_k^*$ discussed in Section \ref{exampleexp} for two values of
$k$ as well as on the optimal weights minimizing the right-hand side of
(\ref{diskret}), where we set $N=100$.
In order to compute the solution $\Xi_{\mathrm{opt}}$, we used the routine
\textit{ipop} from the R-package
\textit{kernlab} by \citet{kernlab2004}. In the
left part of Figure \ref{asymvar}, we show the asymptotic variances of
the different estimators
for the independence
copula. We observe that Pickands estimator has the largest asymptotic variances
(this curve is not displayed in the figure),
while the CFG estimator of \citet{genest2009} yields smaller variances
than the estimator $\hat A_{n,h_1}$, but larger asymptotic variances
than the estimators~$\hat A_{n,h_{5}}$. On the other hand, the estimate
$ \hat A_{n,h_{\mathrm{opt}}}$
corresponding to the numerically determined optimal weight function
yields a substantially smaller variance than
all other estimates under consideration. In the right-hand part
of Figure~\ref{asymvar}, we display the corresponding results
for the asymmetric negative logistic model [see \citet{joe1990}]
%
%
\begin{equation} \label{asyneg}
A(t)=1-\bigl\{ \bigl(\psi_1 (1-t)\bigr)^{-\theta} + (\psi_2t)^{-\theta}
\bigr\}^{-1/\theta}
\end{equation}
with parameters $\psi_1=1,\psi_2=2/3$ and $\theta\in(0,\infty)$
chosen such that the coefficient of tail
dependence is $0.6$.
We observe that the estimate $\hat A_{n,h_5}$ yields the largest
asymptotic variance. The CFG estimate proposed by \citet{genest2009}
and the estimate $\hat A_{n,h_1}$ show a similar behavior (with minor
advantages for the latter), while the best results are
obtained for the new estimate corresponding to the optimal weight function.

We conclude this section with the remark that we have presented a
comparison of the different estimators based on the asymptotic variance which
determines the mean squared error asymptotically.
For finite samples, minimizing only the variance might increase the
bias and therefore the asymptotic
results cannot directly be transferred to applications. In the finite
sample study presented in Section
\ref{simulation}, we will demonstrate that not all of the asymptotic
results yield good predictions for the
finite-sample behavior of the corresponding estimators.

%

\subsection{Convex estimates and endpoint corrections}\label{subsecconvest}

In general, all of the estimates discussed so far [including those
proposed by \citet{genest2009}] will neither be convex, nor will they
satisfy the other characterizing properties of Pickands dependence
functions. However, the literature provides many proposals on how to
enforce these conditions. Various endpoint corrections have been
proposed by \citet{deheuvels1991}, \citet{segers2007}
or \citet{haltaj2000} among others. \citet{filsvilletard2008} proposed
an $L^2$-projection of the estimate
of Pickands dependence function on a space of partially linear
functions which is arbitrarily
close to the space of all convex functions in $\mathcal{A}$ satisfying
the conditions of Lemma \ref{conditionpick}. They also showed that
this transformation
decreases the $L^2$-distance between the ``true'' dependence function
and the estimate.
An alternative concept of constructing convex estimators is based on
the greatest convex
minorant, which yields a decrease in the sup-norm, that is,
\[
\sup_{0<t<1} | \hat A_n^{\mathrm{gcm}} (t) - A(t) | \leq\sup
_{0<t<1} | \hat A_n (t) - A(t) | ,
\]
where $\hat A_n $ is any initial estimate of Pickands dependence
function and
$\hat A_n^{\mathrm{gcm}} $ its greatest convex
minorant [see, e.g., \citet{marshall1970},
\citet{wang1986}, \citet{robertson1988} among others].
It is also possible to combine this concept with an endpoint correction
calculating the greatest convex minorant of the function
\[
t \longrightarrow\bigl(\hat A_n (t) \wedge1\bigr) \vee t \vee(1-t)
\]
[see \citet{genest2009} who also proposed alternative
special endpoint corrections for their estimators]. All these methods
can be used to produce an
estimate of $A$ which has the characterizing properties of a Pickands
dependence function.

\subsection{$M$- and $Z$-estimates} \label{subsecmundzest}

As mentioned in the \hyperref[sec1]{Introduction}, a broader class of
estimates could
be obtained by minimizing more general distances between the given
copula and
the class of functions defined by (\ref{11}) and in this paragraph we briefly
indicate this principle. Consider the best approximation of the copula
$C$ by
functions of the form (\ref{11}) with respect to the distance
%
%
\begin{equation} \label{mabstand}
D_{w}(C,A) = \int_0^1 \int_0^1 \Phi\bigl( C(y^{1-t},y^t),y^{A(t)}
\bigr) w(y,t) \,dy \,dt,
\end{equation}
where $\Phi\dvtx[0,1] \times[0,1] \rightarrow\R^+_0 $ denotes a ``distance''
and $w$ is a given weight function. Note that the minimization in
(\ref{mabstand})
can be carried out by separately minimizing the inner integral for
every value of $t$.
Consequently, the problem reduces to a one-dimensional minimization problem
and assuming differentiability it follows that for fixed $t$ the
optimal value $A^*(t)$
minimizing the interior integral in (\ref{mabstand}) is obtained as a
solution of the equation
\[
\frac{\partial}{\partial a}\int_0^1 \Phi( C(y^{1-t},y^t),y^{a}
) w(y,t) \,dy\bigg|_{a=A^*(t)} =0.
\]
Under suitable assumptions, integration and differentiation can be
exchanged and we have
%
%
\begin{equation} \label{zest}
\int_0^1 \Psi( C(y^{1-t},y^t),y^{a} ) (\log y) y^{a} w(y,t)
\bigg|_{a=A^*(t)} \,dy = 0,
\end{equation}
where $\Psi= \partial_2 \Phi$ denotes the derivative of $\Phi$ with
respect to the second argument.
In general, the solution of (\ref{zest}) is only defined implicitly as
a~functional
of the copula $C$. Therefore, if $C$ is replaced through the empirical
copula the analysis
of the stochastic properties of the corresponding process turns out to
be extremely difficult
because in many cases one has to control improper integrals
(see the proofs of Theorems \ref{theoh1} and \ref{theopickest2} in
the \hyperref[app]{Appendix}).
For the sake of a clear exposition, we do not discuss details in this
paper and defer these
considerations to future research.

Nevertheless, equation (\ref{zest}) yields a different view on the
estimation problem of
Pickands dependence function. Note that the estimate introduced in
Section \ref{subsecnewest} is
obtained by the choice $w(y,t)=h(y)B_h^{-1}$ and
\[
\Phi(z_1,z_2) = (\log z_1 - \log z_2)^2;\qquad\Psi(z_1,z_2) = -2 (\log
z_1 -\log z_2)/z_2
\]
in (\ref{zest}). This estimate corresponds to a minimum distance
estimate. Similarly, an estimate corresponding to the
classical $L^2$-distance is obtained for the
choice
\[
\Phi(z_1,z_2) = (z_1 - z_2)^2;\qquad \Psi(z_1,z_2) = -2 ( z_1 - z_2).
\]
This yields for (\ref{zest}) the equation
\[
\int_0^1 \bigl( C(y^{1-t},y^t) - y^{a} \bigr) (\log y)^2 y^{a} h(
-\log y)
\bigg|_{a=A^*(t)} \,dy = 0,
\]
which cannot be solved analytically. The rank-based versions of
Pickands and the CFG estimator proposed by \citet{genest2009}
do not correspond to $M$-estimates, but could be considered as $Z$-estimates
obtained from (\ref{zest}) for the function
\[
\Psi(z_1,z_2) = (z_1 - z_2)/z_2
\]
with $w_{\mu,\nu} (y) = y^{\mu-1} /(-\log y)^{1+\nu}$ with $\mu
=\nu=0$ and $\mu=0$, $\nu=1$, respectively.
In fact, this choice leads to a general class of estimators
which relates the Pickands and the CFG estimate in an interesting way.
To be precise, note that for $\nu\in[0,1)$ equation (\ref{zest}) yields
%
%
\begin{eqnarray} \label{oneclass}
&&\int_0^1 \frac{( C(y^{1-t},y^t) - I\{ y > e^{-1 }\} ) y^{\mu
-1}}{(-\log y )^\nu} \,dy\nonumber\\
&&\qquad = \int_0^1 \frac{(y^{A(t)} -I\{ y > e^{-1}\}) y^{\mu-1}}{(-\log y
)^\nu} \,dy
\\
&&\qquad = \frac{\Gamma( {1-\nu} )}{(A(t) + \mu)^{1-\nu}} -
\int_0^1 \frac{e^{-\mu x}}{x^\nu} \,dx. \nonumber
\end{eqnarray}
Here the case $\nu=1$ has to be interpreted as the limit $\nu\to1$,
which yields
a~generalization of the defining equation for the CFG estimate, that
is,
\[
- \log\mu- \int_\mu^\infty\frac{e^{-t}}{t}\,dt +
\log\bigl(A(t) + \mu\bigr) =
\int_0^1 \frac{( C(y^{1-t},y^t) - I\{ y > e^{-1 }\} ) y^{\mu
-1}}{\log y } \,dy.
\]
Observing the relation
\[
\lim_{\mu\to0} \log\mu+ \int_\mu^\infty\frac{e^{-t}}{t}\,dt =
-\gamma
\]
we obtain the defining equation for the estimate proposed by \citet
{genest2009} [see
(\ref{cfgdef})]. Similarly, if $\nu\in[0,1)$
it follows from (\ref{oneclass})
%
%
\begin{equation} \label{oneclass1}
\int_0^1 \frac{C(y^{1-t},y^t) y^{\mu-1}}{(-\log y )^\nu} \,dy = \frac
{\Gamma( {1-\nu}) }{(A(t) + \mu)^{1-\nu} }
\end{equation}
and we obtain a defining equation for a generalization of the Pickands
estimate. The classical case
is obtained for $\mu=\nu=0$ [see \citet{genest2009} or equation
(\ref{pickdef})], but (\ref{oneclass1}) defines
many other estimates of this type. Therefore, the Pickands and the CFG
estimate correspond to the extreme
cases in the class $\{ w_{\mu,\nu} |\mu\ge0, \nu\in[0,1] \}$.

We finally note that there are numerous other functions $\Psi$, which
could be used for the construction of
alternative $Z$-estimates, but most of them do not lead to an explicit
solution for
$A^*(t)$. In this sense the CFG-estimator, Pickands-estimator and the
estimates proposed in this paper could be considered
as attractive special cases, which can be explicitly represented in
terms of an integral
of the empirical copula.

\subsection{Finite sample properties} \label{simulation}

In this subsection, we investigate the small sample properties of the new
estimators by means of a simulation study. Especially, we compare the
new estimators with the rank-based estimators suggested by \citet
{genest2009},
which are most similar in spirit with the method proposed in this paper.
We study the finite sample behavior of the
greatest convex minorants of
the endpoint corrected versions of the various estimators.
The new estimators are corrected in a first step by
%
%
\begin{equation} \label{merk}
\hat A_{n,h}^{\mathrm{corr}}(t) := \max\bigl(t,1-t,\min(\hat A_{n,h},1)\bigr)
\end{equation}
and in a second step the greatest convex minorant of $\hat A_{n,h}^{\mathrm{corr}}$
is calculated. For the rank-based CFG and Pickands
estimators, we first used the endpoint corrections proposed in
\citet
{genest2009},
then applied (\ref{merk}) and finally calculated the greatest convex
minorant. Hereby, we compare the
performance of the different statistical procedures which will be used in
concrete applications and
apply the corrections, that are most favorable for the respective
estimators. The
greatest convex minorants are computed using the routine \textit
{gcmlcm} from
the package \textit{fdrtool} by \citet{strimmerconvex}.
All results presented here are based on 5,000 simulation runs and the
sample size is $n=100$.

As estimators, we consider the statistics defined
in (\ref{33}) with the weight function ${h_k}$ and the optimal weight function
determined in Section \ref{subsecoptweight}. An important question is the
choice of the parameter $k$ for the statistic $\hat A_{n,h_k}$ in order
to achieve
a balance between bias and variance. {For this purpose, we first study the
performance of the estimator $\hat A_{n,h_k}$ with respect to different choices
for the parameter $k$ and consider the asymmetric negative logistic model
defined in (\ref{asyneg}) and the symmetric mixed model [see \citet
{tawn1988}]
defined by
%
%
\begin{equation} \label{symmix}
A(t)=1-\theta t+\theta t^2, \qquad\theta\in[0,1].
\end{equation}
The results for other copula models are similar and are omitted for the
sake of brevity.
For the Pickands dependence function (\ref{asyneg}), we used the parameters
$\psi_1=1$ and $\psi_2=2/3$ such that the coefficient of tail
dependence is given by
$\rho=2 (3^\theta+2^\theta)^{-1/\theta}$ and varies in the
interval $(0,2/3)$,
while the parameter $\theta\in[0,1]$ used in (\ref{symmix}) yields
$\rho=\theta/2\in[0,1/2]$.}

\begin{figure}

\includegraphics{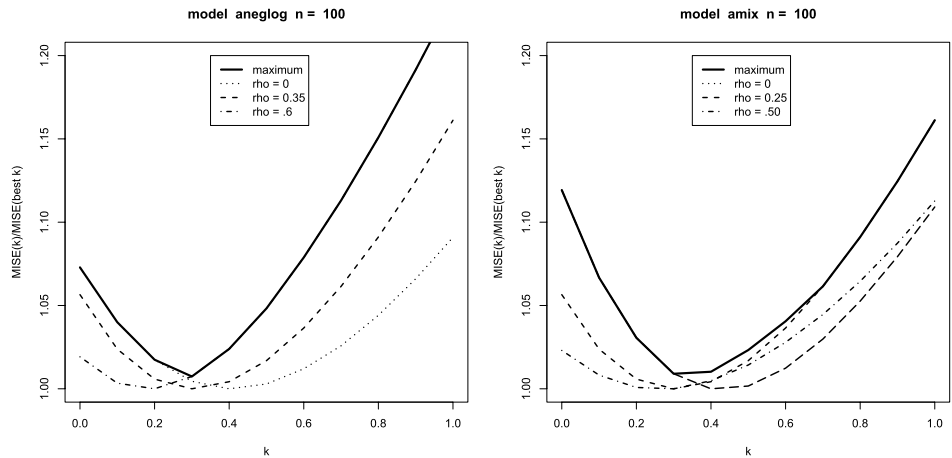}

\caption{The function defined in (\protect\ref{effmise}) for various
models and coefficients of tail dependence. The minimum corresponds to
the optimal value of $k$ in the weight function $h_k$. The solid curve
corresponds to the worst case defined by (\protect\ref{solid}). The sample
size is $n=100$ and the MISE is calculated by 5,000 simulation runs.
Left panel: asymmetric negative logistic model. Right panel: mixed
model.} \label{picMISEk}
\end{figure}

The quality of an estimator $\hat A$ is measured with respect to
mean integrated squared error
\[
\operatorname{MISE} (\hat A) =\mathbb{E}\biggl[ \int_0^1 \bigl(\hat A
(t)-A(t)\bigr)^2 \,dt \biggr],
\]
which was computed by taking the average over 5,000 simulated samples.
The new estimators turned out to be rather robust with respect to the choice
of the parameter $\gamma$ in the definition of the process
$\tilde C_n = \hat C_n \vee n^{-\gamma} $ provided that $\gamma\ge
2/3 $.
For this reason, we use $\gamma=0.95$ throughout this section. Analyzing
the impact of choosing different values for $k$, in Figure \ref{picMISEk}
we display simulated curves
%
%
\begin{equation} \label{effmise}
k \mapsto\frac{\operatorname{MISE}(\hat A_{n,h_k})}{\min_{\ell\ge
0} \operatorname{MISE}(\hat A_{n,h_\ell})}\vadjust{\goodbreak}
\end{equation}
for the asymmetric negative logistic and the mixed models with
different coefficients
of tail dependence $\rho$, as well as the maximum over such curves for
different
values of $\rho$ (solid curves), that is,
%
%
\begin{equation} \label{solid}
k \mapsto\max_\rho\frac{\mathrm{MISE}_\rho(\hat A_{n,h_k})}{\min_{\ell
\geq0}
\mathrm{MISE}_\rho(\hat A_{n,h_\ell})},
\end{equation}
where by $\mathrm{MISE}_\rho$ we denote the MISE for the tail dependence coefficient
$\rho$. The curves in (\ref{effmise}) attain their minima in the optimal
$k$ for the respective~$\rho$, and their shapes provide information
about the
performance of the estimators for nonoptimal values of $k$. The solid
curve gives
an impression about the ``worst case'' scenario (with respect to $\rho
$) in every
model. The simulations indicate, that for $n=100$ the optimal values of
$k$ for
different models and tail dependence coefficients lie in the interval
$[0.2,0.6]$.
Moreover, for values of $k$ in this interval the quality of the
estimators remains
very stable. For $n=200$, $n=500$ and additional models the picture
remains quite similar
and these results are not depicted for the sake of brevity. We thus recommend
using $k = 0.4$ in practical applications. Note that the asymptotic
analysis in
Section \ref{subsecoptweight} suggests that the asymptotically optimal $k$
should differ substantially for various models. However, this effect is
not visible
for sample size up to $n=500$. In these cases, the optimal values for $k$
usually varies in the interval $[0.2,0.8]$.

%

\begin{figure}

\includegraphics{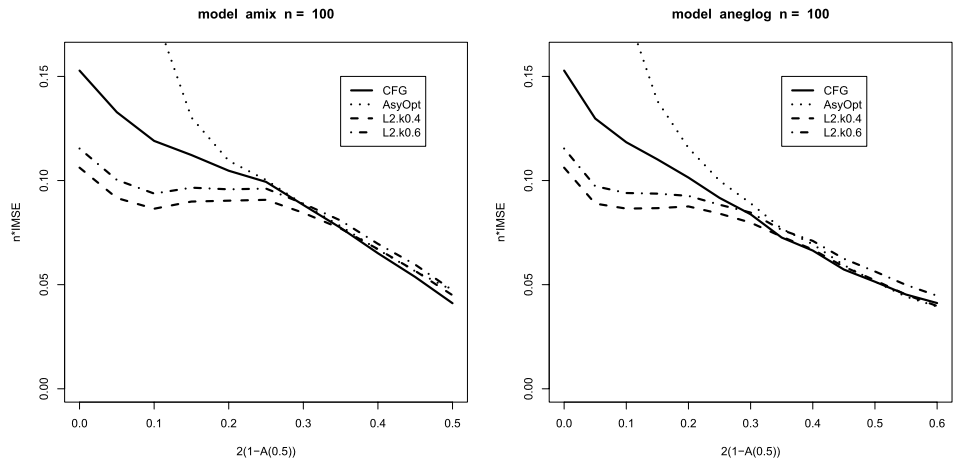}

\caption{$100\times\mathrm{MISE}$ for various estimators, models and
coefficients of tail dependence, based on 5,000 samples of size $n=100$.}
\label{picMISE}
\end{figure}

Next, we compare the new estimators with rank-based versions of Pickands
and the CFG estimator proposed by \citet{genest2009}. In Figure~\ref{picMISE}, the\vadjust{\goodbreak} normalized $\mathrm{MISE}$ is plotted as a
function of the tail
dependence parameter $\rho$ for the asymmetric negative logistic and the
mixed model, where the parameter $\theta$ is chosen in such a way,
that the
coefficient of tail dependence $\rho=2(1-A(0.5))$ varies over the
specific range
of the corresponding model. For each sample, we computed the rank-based
versions of Pickands estimator, the CFG estimator [see \citet{genest2009}]
and two of the new estimators $\hat A_{n,h_k}$ ($k=0.4$, $0.6$).
In this comparison, we also include the estimator $\hat A_{n,h_{\mathrm{opt}}}$ which
uses the optimal weight function determined in Section~\ref{subsecoptweight}.

%
%

Summarizing the results, one can conclude that in general the best
performance is
obtained for our new estimator based on the weight function $h_k$ with $k=0.4$
and $k=0.6$, in particular if the coefficient of tail dependence is
small. A
comparison of the two estimators $\hat A_{n,h_{0.4}}$ and $\hat A_{n,h_{0.6}}$
shows that the choice $k=0.4$ performs slightly better than the choice $k=0.6$
in both models. In both settings, the MISE obtained by $\hat A_{n,h_{0.4}}$
and $\hat A_{n,h_{0.6}}$ is smaller than the MISE of the CFG estimator proposed
in \citet{genest2009} if the coefficient of tail dependence is small.
On the other hand,
the latter estimators yield sightly better results for a large
coefficient of tail dependence.
The results for rank-based version of the Pickands estimator are not depicted,
because this estimator yields a uniformly larger MISE. Simulations of other
scenarios show similar results and are also not displayed for the sake
of brevity.
It is remarkable that the optimal weight function usually yields an
estimator with
a substantially larger MISE than all other estimates if the coefficient
of tail
dependence is small. Similar results can be observed for the sample size
$n=500$ (these results are not depicted). This indicates that the advantages
of the asymptotically optimal weight function only start to play a role
for rather
large sample sizes.\vadjust{\goodbreak}

\section{A test for an extreme-value dependence}\label{sectestevd}

\subsection{The test statistic and its weak convergence}

From the definition of the functional $M_{h}(C,A)$ in (\ref{21}) it is
easy to see that, for a strictly positive weight function $h$ with
$h^*\in L^1(0,1)$,
a copula function $C$ is an extreme-value copula if and only if
\[
\min_{A\in\mathcal{A}}M_{h}(C,A) = M_{h}(C,A^*)=0,
\]
where $A^*$ denotes the best approximation defined in (\ref{astar}).
This suggests to use $M_{h}(\tilde C_{n}, \hat A_{n,h})$
as a test statistic for the hypothesis (\ref{hypextr}), that is,
\[
H_0\dvtx C \mbox{ is an extreme-value copula}.
\]
Recalling the representation (\ref{MCastar})
\[
M_{h}(C,A^*) = \int_0^1\int_0^1 \bar C^2(y,t)h^*(y) \,dy \,dt - B_h
\int_0^1 (A^*(t))^2 \,dt
\]
with $\bar C(y,t)=-\log C(y^{1-t},y^t)$ and defining
$
\bar C_n(y,t) := -\log\tilde C_n(y^{1-t},y^t)
$
we obtain the decomposition
%
%
\begin{eqnarray} \label{decomps1}
&& M_{{h}}(\tilde C_{n}, \hat A_{n,h}) - M_{h}(C,A^*)
\nonumber\\
&&\qquad= \int_0^1\int_0^1 \bigl(\bar C_n^2(y,t) - \bar C^2(y,t)
\bigr)\frac{h^*(y)}{(\log y)^2} \,dy \,dt\nonumber\\
&&\qquad\quad{} - B_h\int_0^1 \hat A_{n,h}^2(t) -
(A^*(t))^2 \,dt
\nonumber\\
&&\qquad= 2\int_0^1\int_0^1 \bigl(\bar C_n(y,t) - \bar C(y,t)
\bigr)\bar C(y,t)\frac{h^*(y)}{(\log y)^2} \,dy \,dt \nonumber\\
&&\qquad\quad{} - 2 B_h\int_0^1 \bigl(\hat A_{n,h}(t) - A^*(t)\bigr) A^*(t) \,dt
\nonumber\\
&&\qquad\quad{} + \int_0^1\int_0^1 \bigl(\bar C_n(y,t) - \bar C(y,t)
\bigr)^2 \frac{h^*(y)}{(\log y)^2} \,dy \,dt\\
&&\qquad\quad{} - B_h\int_0^1 \bigl(\hat A_{n,h}(t)
- A^*(t)\bigr)^2 \,dt
\nonumber\\
&&\qquad= 2\int_0^1\int_0^1 \bigl(\bar C_n(y,t) - \bar C(y,t)
\bigr)\bigl(\bar C(y,t) - A^*(t)(-\log y)\bigr)\frac{h^*(y)}{(\log
y)^2} \,dy \,dt\nonumber\\
&&\qquad\quad{} + \int_0^1\int_0^1 \bigl(\bar C_n(y,t) - \bar C(y,t)
\bigr)^2 \frac{h^*(y)}{(\log y)^2} \,dy \,dt\nonumber\\
&&\qquad\quad{} - B_h\int_0^1 \bigl(\hat A_{n,h}(t) -
A^*(t)\bigr)^2 \,dt\nonumber\\
&&\qquad=: S_1+S_2+S_3, \nonumber
\end{eqnarray}
where the last identity defines the terms $S_1,S_2$ and $S_3$ in an
obvious manner.
Note that under the null hypothesis of extreme-value dependence we have
$A^* = A$ and thus $\bar C(y,t) = A^*(t)(-\log y)$. This means that
under $H_0$ the term $S_1$
will vanish and the asymptotic distribution will be determined by the
large sample properties of the random variable $S_2+S_3$.
Under the alternative, the equality $\bar C(y,t) = A^*(t)(-\log y)$
will not hold anymore and it turns out that in this case the statistic
is asymptotically dominated by the random variable $S_1$.
With the following results, we will derive the limiting distribution of
the proposed test statistic under the null hypothesis and the alternative.
\begin{theorem}\label{testh0}
Assume that the given copula $C$ satisfies condition (\ref{asproc})
and is an extreme-value copula\vspace*{1pt} with Pickands dependence function $A^*$.
If the function $\bar w(y):=h^*(y)/(\log y)^2$ fulfills conditions
(\ref{cw2}) and (\ref{cw3}) for some $\lambda>2$ and
the weight function $h$ is strictly positive and
satisfies assumptions~(\ref{20a}), (\ref{20c}) for $\tilde\lambda:=
\lambda/2>1$, then we have for any
$\gamma\in(1/2, \lambda/4)$ and $n\rightarrow\infty$
\[
n M_{{h}}(\tilde C_n, \hat A_{n,h}) \weak Z_0,
\]
where
the random variable $Z_0$ is defined by
\[
Z_0 := \int_0^1\int_0^1\biggl( \frac{\G
_C(y^{1-t},y^t)}{C(y^{1-t},y^t)}\biggr)^2 \bar w(y) \,dy \,dt - B_h\int
_0^1 \A_{C,h}^2(t) \,dt
\]
with $B_h = \int_0^1 h^*(y)\,dy$ and the process $\{ \A_{C,h}(t)\}
_{t\in[0,1]}$ is defined in Theorem~\ref{theopickest1}.
\end{theorem}

The next theorem gives the distribution of the test statistic
$M_{{h}}(\tilde C_{n}, \hat A_{n,h})$ under the alternative.
Note that in this case we have $M_{h}(C,A^*) >0$.
\begin{theorem}\label{testh1}
Assume that the given copula $C$ satisfies $C\geq\Pi$, condition
(\ref{asproc}) and that $M_{h}(C,A^*) >0$.
If additionally the weight\vspace*{1pt} function $h$ is strictly positive and $h$
and the function $\bar w(y):=h^*(y)/(\log y)^2$ satisfy the assumptions
(\ref{20a}), (\ref{20c})
and (\ref{cw2}), (\ref{cw3})
for some $\lambda>1$, respectively, then
we have for any $\gamma\in(1/2, (1+\lambda)/4\wedge\lambda/2)$ and
$n\rightarrow\infty$
\[
\sqrt n\bigl(M_{h}(\tilde C_{n},\hat A) - M_{h}(C,A^*)\bigr) \weak Z_1,
\]
where
the random variable $Z_1$ is defined as
\[
Z_1 = 2\int_0^1\int_0^1\frac{\G
_C(y^{1-t},y^t)}{C(y^{1-t},y^t)}v(y,t) \,dy \,dt
\]
with
\[
v(y,t) = \bigl( \log C(y^{1-t},y^t) - \log(y)A^*(t) \bigr)\frac{h^*(y)}{(\log y)^2}.
\]
\end{theorem}
\begin{remark}\label{remm}
(a) Note that the weight functions $h_k^*(y) = -y^k \log y$ satisfy the
assumptions of Theorems \ref{testh0} and \ref{testh1} for $k>1$ and $k>0$,
respectively.

(b) The preceding two theorems yield a consistent asymptotic level
$\alpha$ test for the hypothesis of extreme-value dependence by
rejecting the null hypothesis $H_0$ if
%
%
\begin{equation} \label{testext}
n M_{h} (\tilde C_{n},\hat A_{n,h}) > z_{1-\alpha},
\end{equation}
where $z_{1-\alpha}$ denotes the $(1-\alpha)$-quantile of the
distribution of the random variable~$Z_0$.

(c) By Theorem \ref{testh1}, the power of the test (\ref{testext}) is
approximately given by
\begin{eqnarray*}
\PP\bigl( n M_{h} (\tilde C_{n},\hat A_{n,h}) >
z_{1-\alpha}\bigr)
&\approx&1-\Phi\biggl( \frac{z_{1-\alpha}}{\sqrt
{n} \sigma} - \sqrt{n} \frac{M_{h}(C,A^*)}{\sigma} \biggr)\\
&\approx&\Phi\biggl(\sqrt{n} \frac{M_{h}(C,A^*)}{\sigma} \biggr),
\end{eqnarray*}
where the function $A^*$ is defined in (\ref{astar}) corresponding to
the best approximation of the logarithm of the copula $C$ by a function
of the form (\ref{12}),~$\sigma$ is the standard deviation of the distribution of the random variable
$Z_1$ and~$\Phi$ is the standard normal distribution function. Thus,
the power of the test~(\ref{testext})
is an increasing function of the quantity $M_{h}(C,A^*) \sigma^{-1}$.
\end{remark}

\subsection{Multiplier bootstrap} \label{mulitplier}

In general, the distribution of the random variable $Z_0$ cannot be
determined explicitly, because
of its complicated dependence on the (unknown) copula $C$. We hence
propose to determine the quantiles by the multiplier bootstrap approach
as described in
\citet{buecdett2010}.
To be precise, let $\xi_1,\ldots,\xi_n$ denote
independent identically distributed random variables with
\[
\mathbb{P}(\xi_1=0)=\mathbb{P}(\xi_1=2)=1/2.
\]
We define
$\bar\xi_n={n}^{-1}\sum_{i=1}^n \xi_i$ as the mean of $\xi_1,\ldots
,\xi_n$ and consider the
multiplier statistics
\[
\hat C_n^*(u,v) = \hat F_n^*(\hat F_{n1}^{-}(u),\hat F_{n2}^{-}(v)),
\]
where
\[
\hat F_n^*(x_1,x_2) = \frac{1}{n}\sum_{i=1}^n \frac{\xi_i}{\bar\xi
_n} \ind\{X_{i1}\leq x_1, X_{i2}\leq x_2 \},
\]
and $\hat F_{nj}$ denotes the marginal empirical distribution functions.
If we estimate the partial derivatives of the copula $C$ by
\begin{eqnarray*}
\widehat{\partial_1 C}(u,v) &:=& \frac{\hat C_n(u+h,v)-\hat
C_n(u-h,v)}{2h},\\
\widehat{\partial_2 C}(u,v) &:=& \frac{\hat C_n(u,v+h)-\hat C_n(u,v-h)}{2h},
\end{eqnarray*}
where $h=n^{-1/2}\rightarrow0$, we can approximate the distribution of
$\G_C$ by the distribution of the process
%
%
\begin{equation}
\hat\alpha^{\mathrm{pdm}}_n(u,v) := \hat\beta_n(u,v) - \widehat{\partial_1
C}(u,v) \hat\beta_n(u,1)-\widehat{\partial_2 C}(u,v) \hat\beta_n(1,v),\label{remest}
\end{equation}
where $
\hat\beta_n(u,v) =
\sqrt{n} (\hat C_n^*(u,v)-\hat C_n(u,v))$.
More precisely, it was shown by \citet{buecdett2010} that we have weak
convergence conditional on the data in probability toward $\G_C$, that is,
%
%
\begin{equation}\label{weakpdm}
\hat\alpha_n^{\mathrm{pdm}} \weakP\G_C \qquad\mbox{in } l^\infty[0,1]^2,
\end{equation}
where\vspace*{-3pt} the symbol $\displaystyle \weakPi{\xi}$ denotes weak convergence
conditional on the data in probability as defined by \citet
{kosorok2008}, that is, $\alpha_n^{\mathrm{pdm}}\weakPi{\xi}\G_C$ if
%
%
\begin{equation}\label{BL}
\sup_{h \in BL_1(l^\infty[0,1]^2 )} | {\Exp_\xi h(\alpha_n^{\mathrm{pdm}}) -
\Exp h(\G_C)}| \stackrel{\PP}{\longrightarrow} 0
\end{equation}
and
%
%
\begin{equation}\label{am}\quad
\Exp_\xi h(\alpha_n^{\mathrm{pdm}})^*-\Exp_\xi h(\alpha_n^{\mathrm{pdm}})_* \stackrel
{\PP}{\longrightarrow}0 \qquad\mbox{for every } h \in BL_1(l^\infty[0,1]^2).
\end{equation}
Here
\begin{eqnarray*}
&&BL_1(l^\infty[0,1]^2) \\
&&\qquad=\{ f\dvtx l^\infty[0,1]^2
\rightarrow\R \dvtx
\|f\|_\infty\leq1, |f(\beta)-f(\gamma)| \leq\|\beta
-\gamma\|_\infty\\
&&\hspace*{223pt}\forall\gamma,\beta\in l^\infty[0,1]^2 \}
\end{eqnarray*}
is the class of all uniformly bounded functions which are Lipschitz
continuous with constant smaller one,
and $\Exp_\xi$ denotes the conditional expectation with respect to
the weights $\xi_n$ given the data $(X_1,Y_1) \cdots( X_n,Y_n)$.
As a~consequence, we obtain the following bootstrap approximation for $Z_0$.
\begin{theorem} \label{multitest}
If condition (\ref{asproc}) is satisfied, the weight function $h$
satisfies the conditions of
Theorem \ref{testh0} and $h^*(y)(y\log y)^{-2}$ is uniformly bounded then
\begin{eqnarray*}
\hat Z_0^*
&=& \int_0^1\int_0^1\biggl( \frac{\hat\alpha_n^{\mathrm{pdm}}(y^{1-t},y^t)
}{\tilde C_n(y^{1-t},y^t)}\biggr)^2 \frac{h^*(y)}{(\log y)^2} \,dy \,dt
\\
&&{} - B_h^{-1}\int_0^1 \biggl(\int_0^1
\frac{ \hat\alpha_n^{\mathrm{pdm}}(y^{1-t},y^t) }{\tilde C_n
(y^{1-t},y^t)}\frac{h^*(y)}{\log y} \,dy\biggr)^2 \,dt
\end{eqnarray*}
converges weakly to $Z_0$ conditional on the data, that is,
\[
\hat Z_0^*\weakP Z_0 \qquad\mbox{in } l^\infty[0,1].
\]
\end{theorem}

By Theorem \ref{multitest}, $\hat Z_0^*$ is a valid bootstrap
approximation for the distribution of $Z_0$. Consequently, repeating
the procedure $B$ times yields a sample $\hat Z_0^*(1),\ldots,\hat
Z_0^*(B)$ that is approximately distributed according to $Z_0$ and we
can use the empirical $(1-\alpha)$-quantile
of this sample, say
$z_{1-\alpha}^*$, as an approximation for $z_{1-\alpha}$. Therefore,
rejecting the null hypothesis
if
%
%
\begin{equation}\label{testboot}
n M_{h}(\tilde C_n, \hat A_{n,h}) > z_{1-\alpha}^*
\end{equation}
yields a consistent asymptotic level $\alpha$ test for extreme-value
dependence.

Note that the condition on the boundedness of the function
$h^*(y)(y\log y)^2$ is not satisfied for any member of the class
$h^*_k(y)=-y^k / \log(y)$
from Example \ref{exhalphahk}. Nevertheless, mimicking the procedure
from \citet{koyan2010}
and using $h^*_k(y)\ind_{[\varepsilon, 1-\varepsilon]}(y)$ instead
of $h^*_k(y)$ is sufficient for the boundedness. Since this is the procedure
being usually performed in practical applications, Theorem \ref
{multitest} is still valuable for the weight functions investigated in
this paper.

\subsection{Finite sample properties} \label{simualtiontest}
In this subsection, we investigate the finite sample properties of the
test for extreme-value dependence.
We consider the asymmetric negative logistic model (\ref{asyneg}), the
symmetric mixed model~(\ref{symmix})
and additionally the symmetric model of Gumbel
%
%
\begin{equation}
\label{Gumbel}
A(t)=\bigl(t^\theta+(1-t)^\theta\bigr)^{1/\theta}
\end{equation}
with parameter $\theta\in[1,\infty)$ [see \citet{gumbel1960}] and
the model of H\"{u}sler and Reiss
%
%
\begin{equation}\label{huesler}
A(t)=(1-t)\Phi\biggl(\theta+\frac{1}{2 \theta}\log\frac
{1-t}{t}\biggr) + t\Phi\biggl( \theta+\frac{1}{2\theta}\log\frac
{t}{1-t}\biggr),
\end{equation}
where $\theta\in(0,\infty)$ and $\Phi$ is the standard normal
distribution function [see \citet{huesler1989}]. The coefficient of
tail dependence in
(\ref{huesler}) is given by $\rho=2(1-\Phi(\theta))$, that is,
independence is obtained for $\theta\to\infty$ and complete
dependence for $\theta\to0$.
For the Gumbel model (\ref{Gumbel}),
complete dependence is obtained in the limit as $\theta$ approaches
infinity while independence corresponds to
$\theta=1$. The coefficient of tail dependence $\rho=2(1-A(0.5))$ is
given by $\rho=2-2^{1/\theta}$.

We generated 1,000 random samples of sample size $n=200$ from various
copula models and calculated the probability
of rejecting the null hypothesis. Under the null hypothesis, we chose
the model parameters in such a way that the
coefficient of tail dependence $\rho$ varies over the specific range
of the corresponding model. Under the alternative, the coefficient of tail
dependence
does not need to exist and we therefore chose the model parameters,
such that Kendall's $\tau$ is an element of
the set $\{1/4, 1/2, 3/4\}$. The weight function is chosen as
$h_{0.4} (y) = - y^{0.4} / \log(y)$ and the critical values are
determined by the multiplier bootstrap approach as described
in Section \ref{mulitplier} with $B=200$ Bootstrap replications. The
results are stated in Table~\ref{tabtest}.

%
\begin{table}
\caption{Simulated rejection probabilities of
the test (\protect\ref{testboot})
for the null hypothesis of an extreme-value copula for various models.
The first four columns deal with models under the null hypothesis,
while the last four are from the alternative}\label{tabtest}
\begin{tabular}{@{}c@{}}
\begin{tabular}{@{}l@{\hspace*{18pt}}c@{\hspace*{23pt}}c@{\hspace*{23pt}}c@{}}
\hline
$\bolds{H_0}$\textbf{-model}        &  $\bolds{\rho}$ & \textbf{0.05}  & \textbf{0.1}    \\
\hline
Independence       & 0\hphantom{.00}       & 0.031 & 0.075\\[3pt]
Gumbel             & 0.25    & 0.045 & 0.098\\
                   & 0.5\hphantom{0}     & 0.029 & 0.066\\
                   & 0.75    & 0.025 & 0.065\\[3pt]
Mixed model        & 0.25    & 0.043 & 0.09\hphantom{0}\\
                   & 0.5\hphantom{0}    & 0.047 & 0.10\hphantom{0}\\[3pt]
Asy. Neg. Log.     & 0.25    & 0.041 & 0.09\hphantom{0}\\
                   & 0.5\hphantom{0}     & 0.038 & 0.077\\[3pt]
H\"{u}sler--Rei\ss & 0.25    & 0.04\hphantom{0}  & 0.091\\
                   & 0.5\hphantom{0}    & 0.045 & 0.089\\
                   & 0.75    & 0.009 & 0.053\\
\hline
                   \\
\end{tabular}\hspace*{12pt}
\begin{tabular}{@{}l@{\hspace*{18pt}}c@{\hspace*{23pt}}c@{\hspace*{23pt}}c@{}}
\hline
$\bolds{H_1}$\textbf{-model}  & $\bolds{\tau}$    & \textbf{0.05} & \textbf{0.1}\\
\hline
Clayton      & 0.25      & 0.874 & 0.916 \\
             & 0.5\hphantom{0}       & 1\hphantom{.000}     & 1\hphantom{.000} \\
             & 0.75      & 0.999 & 1\hphantom{.000} \\[3pt]
Frank        & 0.25      & 0.291 & 0.396 \\
             & 0.5\hphantom{0}       & 0.73\hphantom{0}  & 0.822 \\
             & 0.75      & 0.783 & 0.898 \\[3pt]
Gaussian     & 0.25      & 0.168 & 0.240 \\
             & 0.5\hphantom{0}       & 0.237 & 0.336 \\
             & 0.75      & 0.084 & 0.156\\[3pt]
$t_4$        & 0.25      & 0.105 & 0.187 \\
             & 0.5\hphantom{0}       & 0.158 & 0.263 \\
             &  0.75     & 0.046 & 0.092 \\
\hline\vspace*{-8pt}
\end{tabular}

\end{tabular}
\end{table}

We observe from the left part of Table \ref{tabtest} that
the level of test is accurately approximated for most of the models, if
the tail dependence is not too strong.
{For a large tail dependence coefficient the bootstrap test is
conservative. This phenomenon can be explained by the fact that for the
limiting case of
random variables distributed according to the upper Fr\'
{e}chet--Hoeffding the empirical copula $\hat C_n$ does not converge
weakly to a nondegenerate process at a rate $1/\sqrt{n}$, rather in
this case it follows that
$\|\hat C_n - C\|=O(1/n)$.
Consequently, the approximations proposed in this paper, which are
based on
the weak convergence of $\sqrt{n} (\hat C_n - C) $ to a nondegenerate
process, are not appropriate for small samples,
if the tail dependence coefficient is large.
Considering the alternative, we observe reasonably good power for the
Frank and Clayton copulas, while for the Gaussian
or $t$-copula deviations from an extreme-value copula are not detected
well with a sample size $n=200$. In some cases, the power
of the test (\ref{testboot}) is close the nominal level. This
observation can be again explained by the closeness to
the upper Fr\'{e}chet--Hoeffding bound.}

Indeed, we can use the minimal distance $M_{h}(C,A^*)$ as a measure
of deviation from an extreme-value copula. Calculating the minimal
distance $M_{h}(C,A^*)$
(with Kendall's $\tau=0.5$ and $h=h_{0.4}$),
we observe that the minimal distances are about ten times smaller for
the Gaussian and $t_4$
than for the Frank and Clayton copula, that is,
%
%
\begin{eqnarray*}
M_{h}(C,A^*_{\mathrm{Clayton}})&=&1.65\times10^{-3},\qquad
M_{h}(C,A^*_{\mathrm{Frank}})= 5.87\times10^{-4},\\
M_{h}(C,A^*_{\mathrm{Gaussian}})&=& 2.08\times10^{-4},\qquad
M_{h}(C,A^*_{t_4})= 1.18\times10^{-4}.
\end{eqnarray*}

Moreover, as explained in Remark \ref{remm}(b) the power of the tests
(\ref{testext}) and~(\ref{testboot}) is
an increasing function of the quantity $p(\mathrm{copula}) =
M_{h}(C,A^*) \sigma^{-1}$. For the four copulas considered in the
simulation study (with $\tau=0.5$), the corresponding ratios are
approximately given by
\begin{eqnarray*}
p(\mathrm{Clayton})&=&0.230,\qquad
p(\mathrm{Frank})=0.134,\\
p(\mathrm{Gaussian})&=&0.083,\qquad
p(t_4)=0.064,
\end{eqnarray*}
which provides some theoretical explanation of the findings presented
in Table \ref{tabtest}.
Loosely speaking, if the value $M_{h}(C,A^*) \sigma^{-1}$ is very
small a larger sample size is required to detect
a deviation from an extreme-value copula. This statement
is confirmed by further simulations results. For example, for the
Gaussian and $t_4$ copula
(with Kendall's $\tau=0.75$) we obtain for the sample size $n=500$ the
rejection probabilities $0.766$ $(0.629)$
and $0.40$ $(0.544)$ for the bootstrap test with level $5\%$ ($10\%$),
respectively.

%

\begin{appendix}\label{app}
\section{Proofs}\label{app1}

\vspace*{-3pt}

\begin{pf*}{Proof of Theorem \ref{theoh1}}
Fix\vspace*{1pt} $\lambda>1$ as in (\ref
{cw3}) and $\gamma\in(1/2,\lambda/2)$.
Due to Lemma 1.10.2(i) in \citet{vandervaart1996}, the process
$\sqrt
n(\tilde C_n - C)$
will have the same weak limit (with respect to the $\weak$
convergence) as $\sqrt n(\hat C_n - C)$.

For $i=2,3,\ldots,$ we consider the following random functions in
$l^\infty[0,1]$:
\begin{eqnarray*}
W_n(t) &=& \int_0^1 \sqrt{n} \bigl( \log\tilde C_n(y^{1-t},y^t) -
\log C(y^{1-t},y^t)\bigr) w(y,t) \,dy, \\
W_{i,n}(t)&=&\int_{1/i}^1 \sqrt{n}\bigl( \log\tilde C_n(y^{1-t},y^t)
- \log C(y^{1-t},y^t)\bigr) w(y,t) \,dy, \\
W(t) &=& \int_{0}^1 \frac{\G_C(y^{1-t}, y^t)}{C(y^{1-t},y^t)} w(y,t)
\,dy, \\
W_i(t) &=& \int_{1/i}^1 \frac{\G_C(y^{1-t}, y^t)}{C(y^{1-t},y^t)}
w(y,t) \,dy.
\end{eqnarray*}
We prove the theorem by an application of Theorem 4.2 in \citet
{billingsley1968}, adapted to the concept of weak
convergence in the sense of Hoffmann--J\o rgensen, see, for example,
\citet{vandervaart1996}. More precisely, we will show in Lemma \ref{lemA1}
in Appendix \ref{seclem} that the weak convergence\vspace*{1pt}
$W_n\weak W$ in $l^\infty[0,1]$ follows from the following three
assertions:
%
%
\begin{eqnarray}\label{cond}
&&\hphantom{\mbox{ii}}\mbox{(i)}\quad \mbox{For every } i\geq2\qquad W_{i,n} \weak W_i
\qquad\mbox{for } n\rightarrow\infty\mbox{ in } l^\infty[0,1], \nonumber\\
&&\hphantom{\mbox{i}}\mbox{(ii)}\quad W_i \weak W\qquad \mbox{for } i\rightarrow\infty
\mbox{ in } l^\infty[0,1],\nonumber\\[-8pt]\\[-8pt]
&&\mbox{(iii)} \quad\mbox{For every }\varepsilon>0\nonumber\\
&&\hphantom{\mbox{(iii)}}\quad\qquad\lim_{i\rightarrow
\infty}\limsup_{n\rightarrow\infty} \PP^*\Bigl({\sup_{t\in[0,1]}}
|W_{i,n}(t)-W_n(t)|>\varepsilon\Bigr) = 0. \nonumber
\end{eqnarray}
The main part of the proof now consists in the verification assertion (iii).

We begin by proving assertion (i). For this purpose, set
$T_i=[1/i,1]^2$ and
consider the mapping
\[
\Phi_1\dvtx\cases{
\mathbb{D}_{\Phi_1} \rightarrow l^\infty(T_i), \cr
f \mapsto\log\circ f,
}
\]
where its domain $\mathbb{D}_{\Phi_1}$ is defined by
$ \mathbb{D}_{\Phi_1}= \{ f\in l^\infty(T_i) \dvtx{\inf_{x\in
T_i}}|f(x)| >0 \} \subset l^\infty(T_i)$.
By Lemma 12.2 in \citet{kosorok2008}, it follows that $\Phi_1$
is Hadamard-differentiable at $C$, tangentially to $l^{\infty}(T_i)$,
with derivative $\Phi'_{1,C}(f) = f/C$.
Since $\tilde C_n \geq n^{-\gamma}$ and $C\geq\Pi$ we have $\tilde
C_n, C\in\mathbb{D}_{\Phi_1}$ and the functional delta method
[see Theorem 2.8 in \citet{kosorok2008}] yields
\[
\sqrt{n}( \log\tilde C_n - \log C) \weak\G_C/C
\]
in $l^\infty(T_i)$. Next, we consider the
operator
\[
\Phi_2\dvtx\cases{
l^\infty(T_i) \rightarrow l^\infty([1/i,1]\times[0,1]),\cr
f \mapsto f\circ\varphi,}
\]
where the mapping
$\varphi\dvtx[1/i,1]\times[0,1]\rightarrow T_i$ is defined by $\varphi
(y,t)=(y^{1-t},y^t)$. Observing
\[
\sup_{(y,t)\in[1/i,1]\times[0,1]} |f\circ\varphi(y,t) - g\circ
\varphi(y,t) | \leq\sup_{\x\in T_i} |f(\x)-g(\x)|
\]
we can conclude that $\Phi_2$ is Lipschitz-continuous. By the
continuous mapping theorem [see, e.g., Theorem 7.7 in \citet{kosorok2008}]
and conditions (\ref{cw1}) and~(\ref{cw2}), we immediately obtain
\[
\sqrt{n}\bigl( \log\tilde C_n(y^{1-t},y^t) - \log C(y^{1-t},y^t)
\bigr) w(y,t) \weak\frac{\G_C(y^{1-t}, y^t)}{C(y^{1-t},y^t)} w(y,t)
\]
in $l^\infty([1/i,1]\times[0,1])$. The assertion in (i) now follows
by continuity of integration with respect to the variable $y$.

For the proof of assertion (ii), we simply note that $\G_C$ is bounded
on $[0,1]^2$ and that
\[
K(y,t)=\frac{w(y,t)}{C(y^{1-t},y^t)}
\]
is
uniformly bounded with respect to $t \in[0,1]$ by the integrable
function $\bar K(y)=\bar w(y) y^{-1}$.

For the proof of assertion (iii), choose some $\alpha\in(0,1/2)$ such
that $\lambda\alpha>\gamma$ and consider the decomposition
%
%
\begin{eqnarray} \label{AiBi}
&&W_n(t)-W_{i,n}(t)\nonumber \\
&&\qquad= \int_0^{1/i} \sqrt{n} \bigl( \log\tilde
C_n(y^{1-t},y^t) - \log C(y^{1-t},y^t)\bigr) w(y,t) \,dy
\\
&&\qquad= B^{(1)}_i(t) + B^{(2)}_i(t), \nonumber
\end{eqnarray}
where
%
%
\begin{equation}\label{Bit}
B^{(j)}_i(t) = \int_{I_{B^{(j)}_i (t)}} \sqrt{n}\log\frac{\tilde
C_n}{C} (y^{1-t},y^t) w(y,t) \,dy,\qquad j=1,2,
\end{equation}
and
%
%
\begin{eqnarray}\label{hd1}
I_{B^{(1)}_i (t)} &=&
\{ 0<y<1/i | C(y^{1-t},y^t)>n^{-\alpha} \},\nonumber\\[-8pt]\\[-8pt]
I_{B^{(2)}_i (t)} &=& (0,1) \setminus I_{B^{(1)}_i (t)}.\nonumber
\end{eqnarray}
The usual estimate
%
%
\begin{eqnarray} \label{iiihalf}
&&\PP^*\Bigl(\sup_{t\in[0,1]} |W_{i,n}(t)-W_n(t)|>
\varepsilon\Bigr)\nonumber\\[-8pt]\\[-8pt]
&&\qquad\leq
\PP^*\Bigl(\sup_{t\in[0,1]}\bigl|B^{(1)}_i(t)\bigr| > \varepsilon
/2\Bigr) + \PP^*\Bigl(\sup_{t\in[0,1]}\bigl|B^{(2)}_i(t)\bigr| > \varepsilon/2\Bigr)\nonumber
\end{eqnarray}
allows for individual investigation of both expressions, and we begin
with the term ${\sup_{t\in[0,1]}} |B^{(1)}_i(t)|$. By the mean value
theorem applied to the logarithm, we have
%
%
\begin{eqnarray} \label{mvt}
\log\frac{\tilde C_n}{C}(y^{1-t},y^t) &=& \log\tilde C_n(y^{1-t},y^t)
- \log{C}(y^{1-t},y^t) \nonumber\\[-8pt]\\[-8pt]
&=& (\tilde C_n - C) (y^{1-t}, {y^t} ) \frac{1}{C^*(y,t)}, \nonumber
\end{eqnarray}
where $C^*(y,t)$ is some intermediate point satisfying $| C^*(y,t)
-C(y^{1-t},y^t) | \le| \tilde C_n(y^{1-t},y^t)- C(y^{1-t},y^t)|$.
Especially, observing $C \geq\Pi$ we have
%
%
\begin{equation}\label{zwischenstelle}
C^*(y,t)\geq(C\wedge\tilde C_n)(y^{1-t},y^t) \geq y\wedge\biggl(y
\frac{\tilde C_n}{C}(y^{1-t},y^t)\biggr)
\end{equation}
and therefore
\begin{eqnarray*}
\sup_{t\in[0,1]}\bigl|B^{(1)}_i(t)\bigr|
&\leq& \sup_{t\in[0,1]} \int_{I_{B^{(1)}_i(t)}} \sqrt{n}|
(\tilde C_n-C) (y^{1-t},y^t) | \\
&&\hspace*{55.5pt}{}\times\biggl| 1\vee\frac
{C}{\tilde C_n} (y^{1-t},y^t)\biggr| w(y,t) y^{-1} \,dy \\
&\leq& \sup_{\x\in[0,1]^2}\sqrt{n}|\tilde C_n(\x)-C(\x)|\\
&&
{}\times
\biggl( 1\vee\sup_{\x\in[0,1]^2 \dvtx C(\x)>n^{-\alpha}} \biggl|
\frac{C}{\tilde C_n}(\x) \biggr| \biggr) \times\psi(i)
\end{eqnarray*}
with $\psi(i)=\int_0^{1/i} \bar w (y) y^{-1} \,dy=o(1)$ for
$i\rightarrow\infty$.
This yields for the first term on the right-hand side of (\ref{iiihalf})
%
%
\begin{eqnarray}
\label{iiisqrt}
\qquad\PP^*\Bigl(\sup_{t\in
[0,1]}\bigl|B^{(1)}_i(t)\bigr| > \varepsilon\Bigr)
& \leq &\PP^*\Biggl(\sup_{\x\in[0,1]^2}\sqrt{n}|\tilde C_n(\x
)-C(\x)| >\sqrt{\frac{\varepsilon}{\psi(i)}} \Biggr)
\nonumber\\[-8pt]\\[-8pt]
\qquad&&{} + \PP^*\Biggl( 1\vee\sup_{C(\x)>n^{-\alpha}} \biggl|
\frac{C}{\tilde C_n}(\x) \biggr| > \sqrt{\frac{\varepsilon}{\psi
(i)}} \Biggr).\nonumber
\end{eqnarray}
Since $\sup_{\x\in[0,1]^2}\sqrt{n}|\tilde C_n(\x)-C(\x)|$ is
asymptotically tight, we immediately obtain
%
%
\begin{equation} \label{510a}
\lim_{i\rightarrow\infty}\limsup_{n\rightarrow\infty} \PP^*
\Biggl(\sup_{\x\in[0,1]^2}\sqrt{n}|\tilde C_n(\x)-C(\x)|>\sqrt{\frac
{\varepsilon}{\psi(i)}} \Biggr) =0.
\end{equation}
For the estimation of the second term in (\ref{iiisqrt}), we
note that
%
%
\begin{equation} \label{Cn-CC}
\qquad\quad \sup_{\x\in[0,1]^2 \dvtx C(\x
)>n^{-\alpha}} \biggl| \frac{\tilde C_n(\x)-C(\x)}{C(\x)} \biggr|
< n^{\alpha} \sup_{\x\in[0,1]^2} |\tilde C_n(\x)-C(\x) |
\stackrel{\PP^*}{\longrightarrow} 0,
\end{equation}
which in turn implies
%
%
\begin{eqnarray} \label{OPbruch1}
\sup_{C(\x)>n^{-\alpha}} \biggl|\frac{C}{\tilde C_n}(\x)\biggr|
&=& \sup_{C(\x)>n^{-\alpha}} \biggl| 1+\frac{\tilde C_n-C}{C}(\x)
\biggr|^{-1} \nonumber\\
&\leq& \biggl( 1-\sup_{C(\x)>n^{-\alpha}} \biggl|\frac{\tilde
C_n-C}{C}(\x) \biggr| \biggr)^{-1} \ind_{A_n} \nonumber\\[-8pt]\\[-8pt]
&&{} + \biggl( \sup_{C(\x)>n^{-\alpha}} \biggl|
1+\frac{\tilde C_n-C}{C}(\x) \biggr|^{-1}\biggr) \ind_{\Omega
\setminus A_n} \nonumber\\
&\stackrel{\PP^*}{\longrightarrow}& 1, \nonumber
\end{eqnarray}
where\vspace*{1pt} $A_n=\{\sup_{C(\x)>n^{-\alpha}} | \frac{\tilde
C_n-C}{C}(\x) | <1/2\}$.
This implies that the function
$\max\{1, {\sup_{C(\x)>n^{-\alpha}}} | \frac{C}{\tilde
C_n}(\x) |\}$
can be bounded by a function that converges to one in outer
probability, and thus
\[
\lim_{i\rightarrow\infty}\limsup_{n\rightarrow\infty} \PP^*
\Biggl( 1\vee\sup_{C(\x)>n^{-\alpha}} \biggl| \frac{C}{\tilde C_n}(\x)
\biggr| > \sqrt{\frac{\varepsilon}{\psi(i)}} \Biggr) = 0.
\]
Observing (\ref{iiisqrt}) and (\ref{510a}) it remains to estimate the
second term on the right-hand side of (\ref{iiihalf}).
We make use of the mean value theorem again [see~(\ref{mvt})]
but use the estimate
%
%
\begin{equation} \label{mvt2}
C^*(y,t)\geq(C\wedge\tilde C_n)(y^{1-t},y^t) \geq y^\lambda\wedge
y^\lambda \frac{\tilde C_n}{C^\lambda}(y^{1-t},y^t)
\end{equation}
[recall that $\lambda>1$ by assumption (\ref{cw3})]. This yields
\begin{eqnarray*}
\sup_{t\in[0,1]}\bigl|B^{(2)}_i(t)\bigr|
&\leq& \sup_{t\in[0,1]} \int
_{I_{B^{(2)}_i(t)}} \sqrt{n}| (\tilde C_n-C) (y^{1-t},y^t)
| \\
&&\hspace*{56.2pt}{}\times\biggl| 1\vee\frac{C^\lambda}{\tilde C_n}
(y^{1-t},y^t)\biggr| w(y,t) y^{-\lambda} \,dy \\
&\leq& \sup_{\x\in[0,1]^2}\sqrt{n}|\tilde C_n(\x)-C(\x)|\\
&&
{}\times
\biggl( 1\vee\sup_{\x\in[0,1]^2 \dvtx C(\x)\leq n^{-\alpha}}
\biggl| \frac{C^\lambda}{\tilde C_n}(\x) \biggr| \biggr) \times\phi(i),
\end{eqnarray*}
where $\phi(i)=\int_0^{1/i} \bar w (y)y^{-\lambda} \,dy =o(1)$ for
$i\rightarrow\infty$ by condition (\ref{cw3}).
Using analogous arguments as for the estimation of ${\sup_{t\in
[0,1]}}|B^{(1)}_i(t)| $ the assertion follows from
\[
\sup_{\x\in[0,1]^2 \dvtx C(\x)\leq n^{-\alpha}} \biggl| \frac
{C^\lambda}{\tilde C_n}(\x) \biggr| \leq\sup_{\x\in[0,1]^2 \dvtx
C(\x)\leq n^{-\alpha}} | n^\gamma C^\lambda(\x) | \leq
n^{\gamma-\lambda\alpha} = o(1)
\]
due to the choice of $\gamma$ and $\alpha$.
\end{pf*}
\begin{pf*}{Proof of Theorem \ref{theopickest2}}
The proof will also be
based on Lemma \ref{lemA1} in Appendix \ref{seclem} verifying
conditions (i)--(iii) in (\ref{cond}).
A careful inspection of the previous proof shows that the verification
of condition (i) in (\ref{cond}) remains valid.
Regarding condition (ii), we have to show that the process $\frac{\G
_C}{C}(y^{1-t},y^t)$ is integrable on the interval $(0,1)$.
For this purpose, we write
\[
\G_C(\x) =\B_C(\x) - \partial_1C(\x)\B_C(x_1,1) - \partial_2
C(\x) \B_C(1,x_2)
\]\looseness=0
and consider each term separately. From Theorem G.1 in \citet
{genest2009}, we know that for any $\omega\in(0,1/2)$ the process
\[
\tilde\B_C(\x) = \cases{
\displaystyle \frac{\B_C(\x)}{(x_1 \wedge x_2)^\omega(1-x_1 \wedge x_2)^\omega},
&\quad if $x_1\wedge x_2\in(0,1)$, \vspace*{2pt}\cr
0, & \quad if $x_1=0$ or $x_2=0$ or
$\x=(1,1)$,}
\]
has continuous sample paths on $[0,1]^2$. Considering
$C(y^{1-t},y^t)\geq y$ and using the notation
%
%
\begin{eqnarray}
\label{K1}
K_1(y,t) &=& {q_\omega(y^{1-t}\wedge y^{t})}y^{-1}, \\
\label{K2}
K_2(y,t) &=& {\partial_1C(y^{1-t},y^t)q_\omega(y^{1-t})}y^{-1}, \\
\label{K3}
K_3(y,t) &=& {\partial_2C(y^{1-t},y^t)q_\omega(y^{t})}y^{-1}
\end{eqnarray}
with $q_\omega(t)=t^\omega(1-t)^\omega$ it remains to show that
there exist integrable functions $K_j^*(y)$ with $K_j(y,t)\leq K_j^*(y)$
for all $t\in[0,1]$ $(j=1,2,3)$.
For $K_1$ this is immediate because $K_1(y,t)\leq(y^{1-t}\wedge
y^{t})^\omega y^{-1}\leq y^{\omega/2-1}$. For $K_2$,\vadjust{\goodbreak} note
that~$\partial_1 C(y^{1-t},y^t)=\mu(t) y^{A(t)-(1-t)}$,
with $\mu(t)=A(t)-tA'(t)$. Therefore,
%
%
\begin{eqnarray} \label{K2*}
\qquad K_2(y,t)\leq\mu(t) y^{A(t)-(1-\omega)(1-t)-1}\leq\mu(t)
y^{\omega/2-1} \leq2 y^{\omega/2-1},
\end{eqnarray}
where the second estimate follows from the inequality $t\vee(1-t)\leq
A(t)\leq1$ and holds for $\omega\in(0,2)$. A similar argument works
for the term $K_3$.

For the verification of condition (iii), we proceed along similar lines
as in the previous proof. We begin by choosing some $\beta\in(1,9/8),
\omega\in(1/4,1/2)$ and some $\alpha\in(4/9,\gamma\wedge(2-\omega
)^{-1})$ in such a\vspace*{1pt} way that $\gamma<\beta\alpha$. First, note that
$y\leq1/(n+2)^2$ implies $\tilde C_n(y^{1-t},y^t)=n^{-\gamma}$ for
all $t\in[0,1]$. This yields
\[
\int_0^{(n+2)^{-2}} \sqrt{n}(\log\tilde C_n-\log C)(y^{1-t},y^t) \,dy
=O\biggl( \frac{\log n}{n^{3/2}} \biggr)
\]
uniformly with respect to $t \in[0,1]$, and therefore it is sufficient
to consider the decomposition in (\ref{AiBi})
with the sets
\begin{eqnarray*}
I_{B^{(1)}_i(t)}&=& \{ 1/(n+2)^{2} <y<1/i |
C(y^{1-t},y^t)>n^{-\alpha} \}, \\
I_{B^{(2)}_i(t)}&=& \bigl(1/(n+2)^{2}, 1/i\bigr) \setminus I_{B^{(1)}_i(t)} .
\end{eqnarray*}
We can estimate the term $B^{(1)}_i(t)$
analogously to the previous proof by
\[
\bigl|B_i^{(1)}(t)\bigr| \leq\int_{I_{B^{(1)}_i(t)}} \sqrt{n}| (\tilde
C_n-C) (y^{1-t},y^t) | \times\biggl| 1\vee\frac{C}{\tilde C_n}
(y^{1-t},y^t)\biggr| y^{-1} \,dy.
\]
Let $H_n$ denote the empirical distribution function of the
standardized sample $ (F(X_{1}), G(Y_1)),\ldots, (F(X_{n}),G(Y_{n}))$.
By the\vspace*{1pt} results in Segers [(\citeyear{segers2010}),
Section~5] we can decompose $
\sqrt{n}(\tilde C_n-C)=\sqrt{n}(C_n\vee n^{-\gamma}-C)$ as
follows:
%
%
\begin{eqnarray}
\label{stute}
\qquad\sqrt{n}(\tilde C_n-C)(\x)
&=& \sqrt{n}(C_n-C)(\x) + \sqrt{n} (\tilde C_n-C_n)(\x)\nonumber
\\
\qquad&=& \alpha_n(\x) - \partial_1C(\x) \alpha_n(x_1,1) - \partial
_2C(\x) \alpha_n(1,x_2) \\
\qquad&&{}+ \tilde R_n(\x),\nonumber
\end{eqnarray}
where
$\alpha_n(\x)=\sqrt{n}(H_n-C)(\x)$ and the remainder satisfies
%
%
\begin{equation} \label{tildeRn}\qquad
{\sup_{\x\in[0,1]^2}} | \tilde R_n(\x) |=O\bigl(n^{1/2-\gamma
}+n^{-1/4}(\log n)^{1/2}(\log\log n)^{3/4}\bigr) \qquad\mbox{a.s.}
\end{equation}
Note that the estimate of (\ref{tildeRn}) requires validity of
condition 5.1 in \citet{segers2010}.
This condition is satisfied provided that the function $A$ is assumed
to be twice continuously differentiable; see Example 6.3 in \citet{segers2010}.
With~(\ref{stute}), we can estimate the term $|B^{(1)}_i(t)| $ analogously
to decomposition (\ref{AiBi}) by $B^{(1)}_{i,1}(t)+\cdots
+B^{(1)}_{i,4}(t)$, where
\begin{eqnarray*}
B^{(1)}_{i,1}(t) &=& \int_{I_{B^{(1)}_i(t)}} | \alpha_n
(y^{1-t},y^t) | \biggl| 1\vee\frac{C}{\tilde C_n}
(y^{1-t},y^t)\biggr| y^{-1} \,dy, \\
B^{(1)}_{i,2}(t) &=& \int_{I_{B^{(1)}_i(t)}} \partial_1C(y^{1-t}, y^t)
| \alpha_n (y^{1-t},1) | \biggl| 1\vee\frac
{C}{\tilde C_n} (y^{1-t},y^t)\biggr| y^{-1} \,dy, \\
B^{(1)}_{i,3}(t) &=& \int_{I_{B^{(1)}_i(t)}} \partial_2C(y^{1-t}, y^t)
| \alpha_n (1,y^t) | \biggl| 1\vee\frac{C}{\tilde
C_n} (y^{1-t},y^t)\biggr| y^{-1} \,dy,\\
B^{(1)}_{i,4}(t) &=& \int_{I_{B^{(1)}_i(t)}} | \tilde R_n
(y^{1-t},y^t) | \biggl| 1\vee\frac{C}{\tilde C_n}
(y^{1-t},y^t)\biggr| y^{-1} \,dy.
\end{eqnarray*}
The decomposition in (\ref{stute}), Theorem G.1 in \citet{genest2009}
and the inequality $\alpha<\gamma\wedge(2-\omega)^{-1}$ may be used
to conclude
\[
\sup_{(y,t)\dvtx C(y^{1-t},y^t)>n^{-\alpha}} \biggl| \frac{\tilde
C_n-C}{C} (y^{1-t},y^t) \biggr| = o_{\PP^*}(1),
\]
which in turn implies
%
%
\begin{equation} \label{OPbruch2}
1\vee\sup_{(y,t)\dvtx C(y^{1-t},y^t)>n^{-\alpha}} \biggl| \frac
{C}{\tilde C_n} (y^{1-t},y^t)\biggr| =O_{\PP^*}(1)
\end{equation}
analogously to (\ref{OPbruch1}). Together with (\ref{tildeRn}) and
observing the inequality $\int_{(n+2)^{-2}}^{1/i} y^{-1} \,dy \leq
2\log(n+2) $, we obtain, for $n\rightarrow\infty$
\[
\sup_{t\in[0,1]} B^{(1)}_{i,4}(t) = O_{\PP^*}\bigl(n^{1/2-\gamma}\log n
+n^{-1/4}(\log n)^{3/2}(\log\log n)^{1/4}\bigr)=o_{\PP^*}(1),
\]
which implies
%
%
\begin{equation} \label{Ai4}
\lim_{i\rightarrow\infty}\limsup_{n\rightarrow\infty} \PP^*
\Bigl( \sup_{t\in[0,1]} B^{(1)}_{i,4}(t) > \varepsilon/4 \Bigr) = 0.
\end{equation}
Observing that $q_\omega(y^{1-t}\wedge y^t) \leq y^{\omega/2}$ the
first term $B^{(1)}_{i,1}(t)$ can be estimated by
\begin{eqnarray*}
\sup_{t\in[0,1]} B^{(1)}_{i,1}(t)
&\leq& \sup_{\x\in[0,1]^2} \frac{|\alpha_n(\x)|}{q_\omega
(x_1\wedge x_2)} \\
&&{}\times
\biggl( 1\vee\sup_{(y,t)\dvtx C(y^{1-t},y^t)>n^{-\alpha}} \biggl|
\frac{C}{\tilde C_n} (y^{1-t},y^t)\biggr| \biggr) \times\psi(i),
\end{eqnarray*}
where $\psi(i)=\int_0^{1/i} y^{-1+\omega/2} \,dy =o(1)$ for
$i\rightarrow\infty$. Using analogous arguments as
in the previous proof we can conclude, using of (\ref{OPbruch2}) and
Theorem G.1 in \citet{genest2009}, that
$ \lim_{i\rightarrow\infty}\limsup_{n\rightarrow\infty} \PP^*(
\sup_{t\in[0,1]} B^{(1)}_{i,1}(t) > \varepsilon/4 ) = 0. $
For the second summand, we note that
\begin{eqnarray*}
\sup_{t\in[0,1]} B^{(1)}_{i,2}(t) &\leq&\sup_{x_1\in[0,1]} \frac
{|\alpha_n(x_1,1)|}{q_\omega(x_1)}
\times
\biggl( 1\vee\sup_{(y,t)\dvtx C(y^{1-t},y^t)>n^{-\alpha}} \biggl|
\frac{C}{\tilde C_n} (y^{1-t},y^t)\biggr| \biggr)\\
&&{} \times\sup_{t\in
[0,1]} \int_0^{1/i} K_2(y,t) \,dy,
\end{eqnarray*}
where $K_2(y,t)$ is defined in (\ref{K2}). Observing the estimate in
(\ref{K2*}), we easily obtain
$ \lim_{i\rightarrow\infty} \sup_{t\in[0,1]} \int_0^{1/i}
K_2(y,t) \,dy =0$.
Again,\vspace*{1pt} under consideration of (\ref{OPbruch2}) and Theorem G.1 in
\citet{genest2009}, we have
$ \lim_{i\rightarrow\infty}\limsup_{n\rightarrow\infty} \PP
^*( \sup_{t\in[0,1]} B^{(1)}_{i,2}(t)\,{>}\,\varepsilon/4 )\,{=}\,0$.
A similar argument works for $B^{(1)}_{i,3}$ and from the estimates for
the different terms the assertion
\[
\lim_{i\rightarrow\infty}\limsup_{n\rightarrow\infty} \PP^*\Bigl(\sup
_{t\in[0,1]}\bigl|B^{(1)}_i(t)\bigr| > \varepsilon\Bigr) =0
\]
follows.
Considering the term ${\sup_{t\in[0,1]} }|B^{(2)}_i(t) |$, we proceed
along similar lines as in the proof of Theorem \ref{theoh1}. For the
sake of brevity, we only state the important differences: in estimation
(\ref{mvt2}) replace $\lambda$ by $\beta$, then make use of
decomposition (\ref{stute}), calculations similar to (\ref{K2*}), and
Theorem G.1 in \citet{genest2009} again and for the estimation of the
remainder note that $\int_{1/(n+2)^2}^{1/i} y^{-\beta}
\,dy=O(n^{2(\beta-1)})$.
\end{pf*}
\begin{pf*}{Proof of Theorem \ref{effic}}
Let $\eta$ denote a
probability measure minimizing the functional $V$
defined in (\ref{asymvarall}). Note that $V$ is convex and define for
$\alpha\in[0,1]$ and a further probability
measure $\xi$ on $[0,1]$ the function
\[
g (\alpha) = V\bigl(\alpha\xi+ (1- \alpha) \eta\bigr).
\]
Because $V$ is convex it follows that $\eta$ is optimal if and only if
the directional derivative of $\eta$ in the direction
$\xi- \eta$ satisfies
\begin{eqnarray*}
0 &\leq& g^\prime(0+) = \lim_{\alpha\to0 +} \frac{g(\alpha)
-g(0)}{\alpha} \\
&=& 2\int^1_0 \int^1_0 k_t (x,y) \,d \xi(x) \,d \eta(y) \\
&&{} - 2 \int^1_0
\int^1_0 k_t (x,y) \,d \eta(x) \,d \eta(y)
\end{eqnarray*}
for all probability measures $\xi$. Using Dirac measures
for $\xi$ yields that this inequality is equivalent to (\ref
{lowbound}), which proves Theorem \ref{effic}.
\end{pf*}
\begin{pf*}{Proof of Theorem \ref{testh0}}
Since the integration mapping is continuous, it suffices to establish
the weak convergence
$
W_n(t) \weak W(t)
$
in $l^\infty[0,1]$ where we define
\begin{eqnarray*}
W_n(t) &=& \int_0^1 n \biggl( \log\frac{\tilde
C_n(y^{1-t},y^t)}{C(y^{1-t},y^t)} \biggr)^2 \bar w(y) \,dy - nB_h\bigl(\hat
A_{n,h}(t) - A^*(t)\bigr)^2,\\
W(t) &=& \int_{0}^1 \biggl(\frac{\G_C(y^{1-t},
y^t)}{C(y^{1-t},y^t)}\biggr)^2 \bar w(y) \,dy - B_h\A_{C,h}^2(t).
\end{eqnarray*}
We prove this assertion along similar lines as in the proof of Theorem
\ref{theoh1}. For $i\geq2$, we recall the notation
$\bar w(y)=h^*(y)/(\log y)^2$ and consider the following random
functions in $l^\infty[0,1]$:
\begin{eqnarray*}
W_{i,n}(t)&=&\int_{1/i}^1 n \biggl( \log\frac{\tilde C_n(y^{1-t},y^t)}
{C(y^{1-t},y^t)}\biggr)^2 \bar w(y) \,dy \\
&&{} - B_h^{-1}\biggl( \int_{1/i}^1 \sqrt n \biggl( \log
\frac{\tilde C_n(y^{1-t},y^t)}{C(y^{1-t},y^t)}\biggr) \frac
{h^*(y)}{\log y} \,dy\biggr)^2,\\
W_i(t) &=& \int_{1/i}^1 \biggl(\frac{\G_C(y^{1-t},
y^t)}{C(y^{1-t},y^t)}\biggr)^2 \bar w(y) \,dy \\
&&{} - B_h^{-1}\biggl(\int_{1/i}^1 \frac{\G_C(y^{1-t},
y^t)}{C(y^{1-t},y^t)} \frac{h^*(y)}{\log y} \,dy\biggr)^2.
\end{eqnarray*}
By an application of Lemma \ref{lemA1} in Appendix \ref{seclem}, it
suffices to show the conditions listed in (\ref{cond}).
By arguments similar to those in the proof of Theorem \ref{theoh1}, we obtain
%
%
\begin{equation}
\sqrt{n}\log\frac{\tilde C_n(y^{1-t},y^t)}{C(y^{1-t},y^t)} \weak
\frac{\G_C(y^{1-t}, y^t)}{C(y^{1-t},y^t)}
\end{equation}
in $l^\infty([1/i,1]\times[0,1])$. Assertion (i) now follows
immediately by the boundedness of the functions $\bar w(y)$ and
$h^*(y)(-\log y)^{-1}$ on $[1/i,1]$
[see conditions~(\ref{cw1}), (\ref{cw2}) and (\ref{20a})] and the
continuous mapping theorem.

For the proof of assertion (ii), we simply note that $\G_C^2$ and $\G
_C$ are bounded on $[0,1]^2$ and
$
K_1(y,t)=\frac{\bar w(y)}{C^2(y^{1-t},y^t)}
$ and $K_2(y,t)=\frac{h^*(y)}{C(y^{1-t},y^t)}$
are bounded\vspace*{1pt} uniformly with respect to $t \in[0,1]$ by the integrable
functions\vspace*{1pt} $\bar K_1(y)=\bar w(y) y^{-2}$ and $\bar K_2(y)=h^*(y)(-\log
y)^{-1} y^{-1}$.

For the proof of assertion (iii), we fix some $\alpha\in(0,1/2)$ such
that $\lambda\alpha>2\gamma$ and consider the decomposition
%
%
\begin{equation} \label{AiBi2}
W_n(t)-W_{i,n}(t) 
= B^{(1)}_i(t) + B^{(2)}_i(t) + B^{(3)}_i(t),
\end{equation}
where
%
%
\begin{eqnarray}
\label{Ait2}
B^{(1)}_i(t) &=& \int_{I_{B^{(1)}_i(t)}} n \biggl( \log\frac{\tilde
C_n(y^{1-t},y^t)}{C(y^{1-t},y^t)} \biggr)^2 \bar w(y) \,dy,
\\
\label{Bit2}
B^{(2)}_i(t) &=& \int_{I_{B^{(2)}_i(t)}} n \biggl( \log\frac{\tilde
C_n(y^{1-t},y^t)}{C(y^{1-t},y^t)} \biggr)^2 \bar w(y) \,dy,
\\
\label{Dit}
B^{(3)}_i(t) &=& -B_h^{-1}I(t,1/i)\bigl(2I(t,1) - I(t,1/i)\bigr),
\end{eqnarray}
$I_{B^{(1)}_i}(t) $ and $I_{B^{(2)}_i}(t) $ are defined in (\ref{hd1}) and
\[
I(t,a) = \sqrt n\int_0^a \biggl( \log\frac{\tilde
C_n(y^{1-t},y^t)}{C(y^{1-t},y^t)} \biggr) \frac{h^*(y)}{\log y} \,dy.
\]
By the same arguments as in the proof of Theorem \ref{theoh1}, we have
for every $\varepsilon>0$
\[
\lim_{i\rightarrow\infty}\limsup_{n\rightarrow\infty} \PP^*\Bigl({\sup
_{t\in[0,1]}} |I(t,1/i)|>\varepsilon\Bigr) = 0,
\]
and
$ \sup_{t\in[0,1]} |I(t,1)| = O_{\PP^*}(1)$,
which yields the asymptotic negligibility\break
$\lim_{i\rightarrow\infty}\limsup_{n\rightarrow\infty} \PP^*(\sup_{t\in[0,1]} |B^{(3)}_i(t)|>\varepsilon) = 0$.
For $B^{(1)}_i(t)$, we obtain the estimate
\begin{eqnarray*}
&&\sup_{t\in[0,1]}\bigl|B^{(1)}_i(t)\bigr|\\
&&\qquad\leq \sup_{t\in[0,1]} \int_{I_{B^{(1)}_i(t)}} n| (\tilde
C_n-C) (y^{1-t},y^t) |^2 \biggl| 1
\vee\frac{C^2}{\tilde C_n^2} (y^{1-t},y^t)\biggr| \bar w(y)
y^{-2} \,dy \\
&&\qquad\leq \sup_{\x\in[0,1]^2}n|\tilde C_n(\x)-C(\x)|^2\times
\biggl( 1\vee\sup_{\x\in[0,1]^2 \dvtx C(\x)>n^{-\alpha}} \biggl| \frac
{C^2}{\tilde C_n^2}(\x) \biggr| \biggr) \times\psi(i),
\end{eqnarray*}
where $\psi(i) := \int_0^{1/i}\bar w(y)y^{-2}\,dy$, which can be
handled by the same arguments as in the proof of Theorem \ref{theoh1}.
Finally, the term $B^{(2)}_i(t)$ can be estimated
by
\begin{eqnarray*}
&& \sup_{t\in[0,1]}\bigl|B^{(2)}_i(t)\bigr| \\
&&\qquad\leq \sup_{t\in[0,1]} \int_{I_{B^{(2)}_i(t)}} n| (\tilde
C_n-C) (y^{1-t},y^t) |^2 \biggl| 1\vee\frac{C^{\lambda
}}{\tilde C_n^2} (y^{1-t},y^t)\biggr| \bar w(y) y^{-\lambda} \,dy
\\
&&\qquad\leq \sup_{\x\in[0,1]^2}n|\tilde C_n(\x)-C(\x)|^2\times
\biggl( 1\vee\sup_{\x\in[0,1]^2 \dvtx C(\x)\leq n^{-\alpha}} \biggl|
\frac{C^{\lambda}}{\tilde C_n^2}(\x) \biggr| \biggr) \times\phi(i),
\end{eqnarray*}
where $\phi(i)=\int_0^{1/i} \bar w(y)y^{-\lambda} \,dy =o(1)$ for
$i\rightarrow\infty$ by condition (\ref{cw3}). Mimic\-king the
arguments from the proof of Theorem \ref{theoh1} completes the proof.~%
\end{pf*}\vspace*{3pt}
\begin{pf*}{Proof of Theorem \ref{testh1}}
Recall the decomposition $M_{h}(\tilde C_{n},\hat A_{n,h}) -
M_{h}(C$, $A^*) = S_1+S_2+S_3$ where $S_1,S_2$ and $S_3$ are
defined in (\ref{decomps1}).
With the notation $\bar v(y):= 2 h^*(y)/(-\log y)$ it follows that
$|v(y,t)|\leq\bar v(y)$ and the assumptions on $h$ yield the validity
of (\ref{cw1})--(\ref{cw3}) for $v(y,t)$. This allows for an
application of Theorem \ref{theoh1} and together with the continuous
mapping theorem we obtain $\sqrt n S_1 \weak Z_1$, where $Z_1$ is the
limiting process defined in (\ref{testh1}). Thus, it remains to verify
the negligibility of $S_2+S_3$. For $S_3$, we note that by Theorem \ref
{theopickest1} and the continuous mapping theorem we have $S_3 = O_{\PP
^*}(1/n)$ and it remains to consider $S_2$. To this end, we fix some
$\alpha\in(0,1/2)$ such that $(1+(\lambda-1)/2)\alpha>\gamma$ and
consider the decomposition\looseness=1
\begin{eqnarray*}
&& \int_0^1 \log^2 \frac{\tilde
C_n(y^{1-t},y^t)}{C(y^{1-t},y^t)}\frac{h^*(y)}{(\log y)^2} \,dy
\\
&&\qquad= \int_{I_{B_1^{(1)}(t)} } \log^2 \frac{\tilde
C_n(y^{1-t},y^t)}{C(y^{1-t},y^t)}\frac{h^*(y)}{(\log y)^2} \,dy\\
&&\qquad\quad{} + \int_{ I_{B_1^{(2)}(t)} } \log^2 \frac{\tilde
C_n(y^{1-t},y^t)}{C(y^{1-t},y^t)}\frac{h^*(y)}{(\log y)^2} \,dy
\\
&&\qquad=: T_1(t,n) + T_2(t,n),
\end{eqnarray*}\looseness=0
where the sets $I_{B_1^{(j)}(t)}, j=1,2$ are defined in (\ref{hd1}).
On the set $I_{B_1^{(1)}(t)}$, we use the estimate
\begin{eqnarray*} 
\log^2 \frac{\tilde C_n(y^{1-t},y^t)}{C(y^{1-t},y^t)} &\leq& \frac
{|\tilde C_n - C|^2}{(C^*)^2}(y^{1-t},y^t) \leq n^\alpha\frac{|\tilde
C_n - C|^2}{C^*}\frac{1}{1\wedge{\tilde C_n}/{C}}(y^{1-t},y^t)
\\
&\leq& n^{\alpha}\frac{|\tilde C_n - C|^2}{C^*}(y^{1-t},y^t)\biggl(1
\vee\sup_{\x\in[0,1]^2 \dvtx C(\x)>n^{-\alpha}}\frac{C(\mathbf
{x})}{\tilde C_n(\mathbf{x})}\biggr),
\end{eqnarray*}
where $| C^*(y,t) -C(y^{1-t},y^t) | \le| \tilde C_n(y^{1-t},y^t)-
C(y^{1-t},y^t)|$.
By arguments similar to those used in the proof of Theorem \ref
{theoh1}, it is now easy to see that
\begin{eqnarray*}
\sqrt n \sup_t |T_1(t,n)| &\leq&\sup_{\x\in[0,1]^2}n^{\alpha
+1/2}|\tilde C_n(\x)-C(\x)|^2\\
&&
{}\times\biggl( 1\vee\sup_{\x\in
[0,1]^2 \dvtx C(\x)>n^{-\alpha}} \biggl| \frac{C}{\tilde C_n}(\x)
\biggr| \biggr)^2 \times K \\
&=& o_{\PP^*}(1),
\end{eqnarray*}
where $K := \int_0^1\bar w(y) y^{-1} \,dy < \infty$ denotes a finite
constant [see condition (\ref{cw3})].
Now set $\beta:= (\lambda-1)/2>0$. From the estimate
\[
C^*(y,t) \geq y^{1+\beta} \biggl( 1 \wedge\frac{\tilde
C_n}{C^{1+\beta}}(y^{1-t},y^t)\biggr) = y^{-\beta} y^\lambda
\biggl(1 \wedge\frac{\tilde C_n}{C^{1+\beta}}(y^{1-t},y^t)\biggr)
\]
we obtain by similar arguments as in the proof of the negligibility of
$|B^{(2)}_i(t)|$ in the proof of Theorem \ref{theoh1} (note that on
$I_{B_1^{(2)}(t)}$\vspace*{-3pt} we have $y \leq C(y^{1-t},\break y^t) \leq n^{-\alpha}$ )
\begin{eqnarray*}
{\sup_{t \in[0,1]}} |T_2(t,n)| &\leq&\log(n) n^{-\beta\alpha}\sup
_{\x\in[0,1]^2}\sqrt{n}|\tilde C_n(\x)-C(\x)|\\
&&{}\times\biggl( 1\vee\sup_{\x\in[0,1]^2 \dvtx
C(\x)\leq n^{-\alpha}} \biggl| \frac{C^{1+\beta}}{\tilde C_n}(\x)
\biggr| \biggr) \times\tilde K,
\end{eqnarray*}
where\vspace*{1pt} $\tilde K := \gamma\int_0^1 (1-\log y ) \frac{h^*(y)}{(\log
y)^2} y^{-\lambda} \,dy $ denotes a finite constant [see conditions
(\ref{cw3}) and (\ref{20c})] and we used the estimate
\begin{eqnarray*}
\biggl|{\log\frac{\tilde C_n(y^{1-t},y^t)}{C(y^{1-t},y^t)}}\biggr|^2
&\leq&(\gamma\log n - \log y )\biggl|{\log\frac{\tilde
C_n(y^{1-t},y^t)}{C(y^{1-t},y^t)}}\biggr| \\
&\leq&\gamma\log(n) (1 - \log y )\biggl|{\log\frac{\tilde
C_n(y^{1-t},y^t)}{C(y^{1-t},y^t)}}\biggr|,
\end{eqnarray*}
which holds for sufficiently large $n$.
Finally, we observe that
\[
\sup_{\x\in[0,1]^2 \dvtx C(\x)\leq n^{-\alpha}} \biggl| \frac
{C^{1+\beta}}{\tilde C_n}(\x) \biggr| \leq\sup_{\x\dvtx C(\x)\leq
n^{-\alpha}} | n^\gamma C^{1+\beta}(\x) | \leq
n^{\gamma-(1+\beta)\alpha} = o(1).
\]
Now the proof is complete.
\end{pf*}
\begin{pf*}{Proof of Theorem \ref{multitest}}
The conditions on the weight function imply that all integrals in the
definition of $Z_0$ are proper and therefore the mapping $(\G
_C,C)\mapsto Z_0(\G_C,C)$ is continuous. Hence, the result follows by
the continuous mapping theorem for the bootstrap [see, e.g., Theorem
10.8 in \citet{kosorok2008}] provided the conditional weak convergence
in (\ref{weakpdm}) holds under the nonrestrictive smoothness
assumption (\ref{asproc}). To see this, proceed similar as in
\citet
{buecdett2010} and show Hadamard-differentiability of the mapping
$H\mapsto H(H_1^-,H_2^-)$, which is defined for some distribution
function $H$ on the unit square whose marginals $H_1=H(\cdot,1)$ and
$H_2=H(1,\cdot)$ satisfy $H_1(0)=H_2(0)=0$. This can be done by
similar arguments as in \citet{segers2010} and the details are omitted
for the sake of brevity.
\end{pf*}

\section{An auxiliary result}\label{seclem}

\begin{lemma} \label{lemA1}
Let $X_n, X_{i,n}\dvtx\Omega\rightarrow\Dm$ for $i,n\in\Nat$ be
arbitrary maps with values in the metric space $(\Dm,d)$ and
$X_i,X\dvtx\Omega\rightarrow\Dm$ be Borel-measurable. Suppose that:
\begin{eqnarray*}
&&\hphantom{\mbox{\textup{ii}}}\mbox{\textup{(i)}}\quad  \mbox{For every } i\in\Nat\qquad X_{i,n} \weak X_i
\qquad\mbox{for } n\rightarrow\infty, \\
&&\hphantom{\mbox{\textup{i}}}\mbox{\textup{(ii)}}\quad X_i \weak X \qquad\mbox{for } i\rightarrow\infty, \\
&&\mbox{\textup{(iii)}}\quad \mbox{For every }\varepsilon>0\qquad \lim_{i\rightarrow
\infty}\limsup_{n\rightarrow\infty} \PP
^*\bigl(d(X_{i,n},X_n)>\varepsilon\bigr) = 0.
\end{eqnarray*}
Then $X_n\weak X$ for $n\rightarrow\infty$.\vadjust{\goodbreak}
\end{lemma}
\begin{pf}
Let $F\subset\Dm$ be closed and fix $\varepsilon
>0$. If $F^\varepsilon=\{x\in\Dm\dvtx d(x,F)\leq\varepsilon)$ denotes
the $\varepsilon$-enlargement of $F$ we obtain
\[
\PP^*(X_n\in F)\leq\PP^*(X_{i,n}\in F^{\varepsilon}) + \PP
^*\bigl(d(X_{i,n},X_n)>\varepsilon\bigr).
\]
By hypothesis (i) and the Portmanteau theorem [see \citet{vandervaart1996}]
\[
\limsup_{n\rightarrow\infty}\PP^*(X_n\in F)\leq\PP(X_{i}\in
F^{\varepsilon}) +\limsup_{n\rightarrow\infty}\PP
^*\bigl(d(X_{i,n},X_n)>\varepsilon\bigr).
\]
By conditions (ii) and (iii)
$\limsup_{n\rightarrow\infty} \PP^*(X_n\in F) \leq P(X\in
F^\varepsilon)$
and since $ F^\varepsilon\downarrow F$ for $\varepsilon\downarrow0$
and closed $F$ the result follows by the Portmanteau theorem.
\end{pf}
\end{appendix}

\section*{Acknowledgments}
The authors would like to thank Martina Stein, who typed parts of this
manuscript with considerable technical expertise.
The authors would also like to thank Christian Genest for pointing out
important references and Johan Segers
for many fruitful discussions on the subject.
We are also grateful to two unknown referees and an Associate Editor
for their constructive comments on an earlier version of this
manuscript, which led to a substantial improvement of the paper.


%

%
\printaddresses

\end{document}